%% file: main.tex
\documentclass[a4paper,12pt]{article}

\usepackage[english]{babel}
\usepackage[utf8x]{inputenc}
\usepackage[T1]{fontenc}

\usepackage[a4paper,top=2cm,bottom=2cm,left=1.75cm,right=1.75cm,marginparwidth=1.75cm]{geometry}
\usepackage{subfiles}
\usepackage{amsmath}
\usepackage{commath}
\usepackage{calrsfs}
\usepackage{euscript}
\usepackage{standalone}
\usepackage{enumerate}
\usepackage{graphicx}
\usepackage{color}
\usepackage[numbers]{natbib}
\usepackage{appendix}
\usepackage{graphicx}
\usepackage{pgfplots}
\usepackage{tikz}
\usetikzlibrary{shapes.geometric, arrows}
\pgfplotsset{compat=1.14}
\usepackage{xfrac}
\usepackage{subcaption}
\usepackage{array}
\usepackage{booktabs}
\usepackage{authblk}
\usepackage{lipsum}

\title{A Preconditioned Multiple Shooting Shadowing Algorithm for the Sensitivity Analysis of Chaotic Systems}
\author{Karim Shawki\footnote{Email: karim.shawki14@imperial.ac.uk} }
\author{George Papadakis\footnote{Email: g.papadakis@imperial.ac.uk}}
\affil{Department of Aeronautics, Imperial College London, Exhibition Road, London, SW7 2AZ, UK}

\begin{document}
\maketitle
\begin{abstract}
We propose a preconditioner that can accelerate the rate of convergence of the Multiple Shooting Shadowing (MSS) method \cite{Blonigan2018}. This recently proposed method can be used to compute derivatives of time-averaged objectives (also known as sensitivities) to system parameter(s) for chaotic systems. We propose a block diagonal preconditioner, which is based on a partial singular value decomposition of the MSS constraint matrix. The preconditioner can be computed using matrix-vector products only (i.e. it is matrix-free) and is fully parallelised in the time domain. Two chaotic systems are considered, the Lorenz system and the 1D Kuramoto Sivashinsky equation. Combination of the preconditioner with a regularisation method leads to tight bracketing of the eigenvalues to a narrow range. This combination results in a significant reduction in the number of iterations, and renders the convergence rate almost independent of the number of degrees of freedom of the system, and the length of the trajectory that is used to compute the time-averaged objective. This can potentially allow the method to be used for large chaotic systems (such as turbulent flows) and optimal control applications. The singular value decomposition of the constraint matrix can also be used to quantify the effect of the system condition on the accuracy of the sensitivities. In fact, neglecting the singular modes affected by noise, we recover the correct values of sensitivity that match closely with those obtained with finite differences for the Kuramoto Sivashinsky equation in the light turbulent regime.
\makeatletter{\renewcommand*{\@makefnmark}{}
\footnotetext{Declarations of interest: none}\makeatother}
\end{abstract}
\section{Introduction}
\input{Introduction}

\section{Shadowing Methods}\label{Sec1}
\input{Sec1}

\section{A Preconditioner for the MSS Schur Complement}\label{Schur-PC}
\input{Schur-PC}

\section{Application to the Lorenz System}\label{PC-Lorenz}
\input{PC-Lorenz}
\clearpage

\section{Application to the Kuramoto Sivashinsky Equation}\label{BDP_KS_App}
\input{BDP_KS_App}

\section{Effect of the System Condition on the Accuracy of the Computed Sensitivity}\label{ill_cond}
\input{ill_cond}

\section{Regularization of the Preconditioned System}\label{Tikhonov}
\input{Tikhonov}

\section{Computational Cost of the Method}\label{cost}
\input{cost}

\section{Conclusions}\label{conclusions}
\input{Conclusions}

\section*{Acknowledgements}
The first author wishes to acknowledge the financial support of Al-Alfi Foundation in the form of a PhD scholarship. 
\clearpage

\bibliographystyle{unsrt}
\bibliography{main.bib}
\end{document}

%% file: Introduction.tex
Optimisation has a broad range of applications. In aerospace engineering for example, optimisation can be used in airfoil design \cite{Jameson1988,Liao2010}, flow control \cite{Xiao2017NonlinearFlow,Luchini2014AdjointAnalysis,Peitz2015} and uncertainty quantification \cite{Liao2013UncertaintyNonlinearities}. In practice, we are usually interested in optimising a time-averaged quantity (objective) $\bar{J}$ subject to a set of constraints. An example is the minimisation of an airfoil drag by some active control technique, subject to the governing equations, boundary and initial conditions. The sensitivity $\sfrac{d\bar{J}}{ds}$, where $s$ is one or more control parameters, can be employed together with a gradient-based algorithm to find the optimal values of $s$ that minimize the objective $\bar{J}$.

Traditional sensitivity analysis methods involve solving the tangent equation (or the adjoint equations for multiple control parameters $s$) to compute $\sfrac{d\bar{J}}{ds}$. The above methods are based on linearisation around a reference trajectory in phase space.  If the dynamical system is chaotic, $\sfrac{d\bar{J}}{ds}$ diverges exponentially for long time horizons $T$. Here `long time horizons' refers to time intervals much longer than the characteristic time scales of the system dynamics. Lea et al. \cite{Lea2002SensitivitySystem} showed that the adjoint of the Lorenz system yielded $|\sfrac{d \bar{J}}{d s}|\propto e^{\lambda_{max} t}$, where $\lambda_{max}$ is the largest Lyapunov exponent, that describes the dominant rate of separation of two trajectories evaluated with parameters $s$ and $s+\delta s$. Transitional and turbulent flows are chaotic; therefore, traditional methods fail to compute the correct derivatives for long $T$. A recent example is the non-linear optimal control of bypass transition in a boundary layer using blowing and suction, where $T$ had to be restricted to avoid the exponential growth of the adjoint variables during backward integration \cite{Xiao2017NonlinearFlow}.

Many attempts have been made to develop methods that can compute accurate sensitivities $\sfrac{d\bar{J}}{ds}$ for long $T$. In \cite{Lea2002SensitivitySystem}, the `ensemble-adjoint' method was proposed, where a long $T$ is split into many integration segments. For a given $s$, the adjoint method is applied to all segments, and an average is taken. The accuracy of  $\sfrac{d\bar{J}}{ds}$ depends on the length of segments, and knowledge of the suitable lengths is not known a priori. 

Thuburn \cite{Thuburn2005} employed the Fokker-Planck equation, which governs the evolution of the probability density function in phase space when stochastic Wiener forcing terms are added to the discrete model equations. These terms are artificial but they introduce diffusion to the evolution equation, thus ensuring a smooth steady solution. He expressed the long-time sensitivity in terms of the adjoint of the Fokker-Plank equation. The issue with this method however, is its computational cost for large systems (a multi-dimensional diffusion operator must be inverted) and solution accuracy (due to the addition of the stochastic terms).

The Fluctuation-Dissipation Theorem, a powerful tool of statistical physics \cite{Kubo_1966}, can also be used to estimate the time-average response to forcing in nonlinear systems. The key underlying assumption is that the forcing is weak enough such that the response of the system changes linearly with the forcing. A 1D reduced order model (based on Dynamic Mode Decomposition) was recently proposed to estimate the response to external forcing of the horizontally averaged temperature in a fully turbulent, buoyancy driven flow  \cite{Khodkar2018}. The model was also used to derive the forcing that elicits a desired response, and this makes it suitable for control applications.

Lasagna \cite{Lasagna_2018} formulated a well-behaved sensitivity analysis method for unstable periodic orbits. Due to periodicity, the adjoints are bounded in time irrespective of the orbit length, and the sensitivities computed for every orbit are exact. However, finding periodic orbits can be very challenging for 3D turbulent flows.

The Least Squares Shadowing (LSS) method \cite{Wang2014} considers a reference trajectory in phase space, evaluated at $s$, and finds a nearby (or `shadowing trajectory'), evaluated at $s+\delta s$, that stays close to the reference trajectory. The proximity of the two trajectories regularises the problem and can be used to find accurate $\sfrac{d\bar{J}}{ds}$ for long $T$. The method relies on the shadowing lemma, that guarantees that such a nearby trajectory exists for structurally stable, uniformly hyperbolic systems \cite{Bowen1975-LimitDiffeomorphisms,Pilyugin1999ShadowingSystems,Holmes2012TurbulenceSymmetry}. The two trajectories start from different initial conditions, but for ergodic systems, the time-average and its sensitivity do not depend on the initial conditions. LSS results in a linear two-point boundary value problem in time, which is time-consuming to solve and has large storage requirements.

Blonigan and Wang \cite{Blonigan2018} proposed the Multiple Shooting Shadowing (MSS) algorithm that results in a much smaller linear system than that of LSS. This is because the norm of the distance between the reference and shadowing trajectories is defined using the values at  discrete checkpoints, not at every time step. The system can be solved iteratively using Krylov subspace methods (such as Conjugate Gradient) that only require matrix-vector products. Although computationally more efficient than standard LSS, in practice the convergence rate can be slow due to the high condition number of the system matrix.

Non-Intrusive Least Squares Shadowing (NILSS) \cite{Ni2017} was proposed to further reduce the computational cost, and it is closely related to MSS \cite{Blonigan2018}. NILSS requires the integration of one inhomogeneous tangent that represents the effect of changing $s$, and a set of homogeneous tangents that represent the effect of changing initial conditions. The reference trajectory is segmented and a minimisation problem is solved for the set of weights of the homogeneous tangent solutions that cancel out the growth of the inhomogeneous solution at every segment, leading to an overall bounded result. The cost of solving the NILSS matrix is negligible compared to the cost of solving the Transcription LSS matrix. The main cost is in setting up the matrix, which requires many integrations of the tangent equation, at least equal to the number of positive Lyapunov exponents. For high $Re$ flows the cost of tangent equation integrations would be high due to the expected large number of positive Lyapunov exponents. NILSS has been applied successfully to a chaotic flow over a backwards facing step \cite{Ni2017}. A discrete adjoint version was derived in \cite{Blonigan2017} and applied for the first time to a 3D DNS simulation of turbulent channel flow at $Re_{\tau}=140$.  

MSS remains an attractive algorithm because it does not require the prior computation of all the positive Lyapunov exponents. However, it suffers from slow convergence, as already mentioned. In this paper, we propose a preconditioner that can significantly accelerate the convergence rate. We show that the operations to construct the preconditioner are fully parallel-in-time, require only matrix-vector products, and lead to significant overall cost savings (compared to standard MSS). Using this preconditioner and regularising the small eigenvalues of the system matrix, leads to accurate results and convergence rates almost independent of $T$, and most importantly, of the number of degrees of freedom, $N$. 

This paper is structured as follows: Section (\ref{Sec1}) presents a general overview of LSS and MSS. In Section (\ref{Schur-PC}), the new preconditioner for MSS is derived and is applied to the Lorenz system in Section (\ref{PC-Lorenz}), and the Kuramoto Sivashinsky equation in Section (\ref{BDP_KS_App}). The effect of ill-conditioning on the solution accuracy is explored in Section (\ref{ill_cond}) and the regularisation of the preconditioned system is investigated in Section (\ref{Tikhonov}). The computational cost of the proposed method is discussed in Section (\ref{cost}) and we conclude in section (\ref{conclusions}).

%% file: Sec1.tex
Consider a dynamical system described by the set of ordinary differential equations
\begin{equation}\label{ODE eq}
\frac{d\textbf{u}}{dt}=\textbf{f}(\textbf{u},s), \quad u(0,s)=u_0(s)
\end{equation}
where $\textbf{u}(t,s)$ is the vector of state variables of length $N$, $\textbf{f}$  is the vector of non-linear equations, and $s$ is one (or more) system or control parameters. We define a time-averaged objective as
\begin{equation}
\label{J avg eq}
\bar{J}=\frac{1}{T}\int_{0}^{T} J(\textbf{u}(t,s),s) \, dt
\end{equation}
and we are interested in computing the sensitivity $\sfrac{d\bar{J}}{ds}$ over a long time horizon, $T\rightarrow \infty$. Using the chain rule, we can write
\begin{equation}\label{sensitivity eqn}
\frac{d\bar{J}}{ds}=\frac{1}{T}\int_{0}^{T}  \frac{\partial J(\textbf{u}(t,s),s)}{\partial u}\textbf{v}(t,s) +\frac{\partial J}{\partial s}\, dt
\end{equation}
where $\textbf{v}(t,s)=\sfrac{d u}{d s}$. Perturbing  \eqref{ODE eq} by $\delta s$ and linearizing, the tangent evolution equation for $\textbf{v}(t,s)$ can be obtained:
\begin{equation}\label{tangent eqn}
\frac{d\textbf{v}}{dt}=\frac{\partial f}{\partial u}\textbf{v}+\frac{\partial\textbf{f}}{\partial s}
\end{equation}
with the initial condition $\textbf{v}(0)=\sfrac{du_0}{ds}$. Equation \eqref{tangent eqn} is solved for $\textbf{v}(t)$ and the sensitivity is found using \eqref{sensitivity eqn}. An adjoint version of \eqref{tangent eqn} can be derived to compute the derivative of $\bar{J}$ to many parameters simultaneously.

Chaotic systems of equations are sensitive to small changes in initial conditions and parameter values. For such systems, the solution $\textbf{v}(t)$ grows exponentially at a rate dictated by $\lambda_{max}$, the maximum positive Lyapunov exponent. Since all chaotic systems have at least one positive $\lambda$, the sensitivity computed through \eqref{sensitivity eqn} becomes meaningless for long $T$. Shadowing methods ensure that trajectories evaluated at $s$ and $s+\delta s$ do not diverge, and therefore the solution to \eqref{tangent eqn} remains bounded in time.

\subsection{Least Squares Shadowing (LSS)}
The Least Squares Shadowing method (LSS) was developed by Wang et al. \cite{Wang2014} to compute the sensitivity $\sfrac{d\bar{J}}{ds}$ over long $T$ for chaotic systems. LSS finds a `shadow' trajectory that satisfies (\ref{ODE eq}) at a slightly perturbed parameter value, i.e. $\sfrac{d\textbf{u}'}{dt}=\textbf{f}(\textbf{u}',s+\delta s)$, that stays close in phase space  to a reference trajectory $\textbf{u}_{ref}$ evaluated at $s$, for all $T$. Since the solution  $\textbf{u}'(\tau(t))$ does not diverge from $\textbf{u}_{ref}(t)$, it can be used to evaluate meaningful sensitivities $\sfrac{d\bar{J}}{ds}$. A time transformation $\tau(t)$ is required to avoid algebraically growing components arising from the projection of $\textbf{v}(t)$ on the direction of $\textbf{f}(\textbf{u}(t))$ \cite{Wang2013}. To obtain the shadowing trajectory, a minimisation problem is formulated and solved. Linearisation results in
\begin{equation} \label{LSS lin eq}
\begin{aligned}
\underset{\textbf{v},\eta}{\text{Minimise}} \, \, \,  \,      & \frac{1}{T}\int_{0}^{T} \|\textbf{v}\|^2 + \alpha^2\eta^2  \, dt, \\
\text{subject} \, \text{to} \, \, \, \,  & \frac{d\textbf{v}}{dt}=\frac{\partial f}{\partial u}\textbf{v}+\frac{\partial \textbf{f}}{\partial s}+\eta \textbf{f}(u_r,s) 
\end{aligned}
\end{equation}

\noindent where  $\eta(t)=\sfrac{d\left(\sfrac{d\tau}{dt}-1\right)}{ds}$ is a time dilatation term, and $\textbf{v}(t)=\sfrac{d(\textbf{u}'-\textbf{u}_{ref})}{d s}$. The solutions $\textbf{v}(t)$ and $\eta(t)$ are then used to compute the sensitivity by evaluating (see \cite{Wang2014} for the derivation):
\begin{equation}\label{dJ/ds LSS eq}
\frac{d \bar{J}}{d s}=\frac{1}{T}\int_{0}^{T}\left(\frac{\partial J}{\partial u}\textbf{v}+\frac{\partial J}{\partial s}+\eta(J-\bar{J})\right) \, dt
\end{equation}
The solution of \eqref{LSS lin eq} results in a two-point boundary value problem in time. Discretisation in space and time yields a linear system with $M\times N$ unknowns, where $M$ is the number of time steps, and $N$ is the number of degrees of freedom. The large size of the matrix therefore makes the method applicable to low-dimensional systems only.

\subsection{Multiple Shooting Shadowing (MSS)}
In attempt to reduce the cost of LSS for high dimensional systems, Blonigan and Wang \cite{Blonigan2018} proposed to minimise $\textbf{v}(t)$ at $K$ discreet checkpoints as shown in Figure (\ref{multiseg diag}). This idea reduces the number of unknowns to $K \times N \ll M\times N$. The Multiple Shooting Shadowing (MSS) method is formulated as:
\begin{subequations}\label{MSS eqn}
\begin{align}
\underset{\textbf{v}(t_i^+)}{\text{Minimise}}  &\, \, \, \sum_{i=0}^K \| \textbf{v}(t_i^+) \|_2^2 \\
\text{subject} \, \text{to} & \, \, \, \textbf{v}(t_i^+)=\textbf{v}(t_i^-)  &(i=1,2,...,K-1) \\
&\, \, \, \frac{d\textbf{v}}{dt}-\frac{\partial f}{\partial u}\textbf{v}-\frac{\partial \textbf{f}}{\partial s}-\eta \textbf{f} =0 & t_i<t<t_{i+1}  \, \, \, (i=0,1,...,K-1)\\
&\langle \textbf{f}(\textbf{u}(t),s),\textbf{v}(t) \rangle=0 & t_i<t<t_{i+1} \,\,\, (i=0,1,...,K-1)
\end{align}
\end{subequations}
The constraint (\ref{MSS eqn}b) enforces the continuity of $\textbf{v}(t)$ between the inner checkpoints ($t_1,t_2,...,t_{K-1}$). The inner product (\ref{MSS eqn}d) ensures that $\textbf{v}(t)$ remains normal to $\textbf{f}(\textbf{u}(t),s)$. This is required because components of $\textbf{v}(t)$ not normal to $\textbf{f}(\textbf{u}(t),s)$ lead to linear growth in $\textbf{v}(t)$ \cite{Blonigan2017}.

\tikzstyle{icblock} = [rectangle, rounded corners, minimum width=2cm, minimum height=1cm,text centered, draw=black]
\tikzstyle{tblock} = [rectangle, rounded corners, minimum width=2cm, minimum height=1cm,text centered]
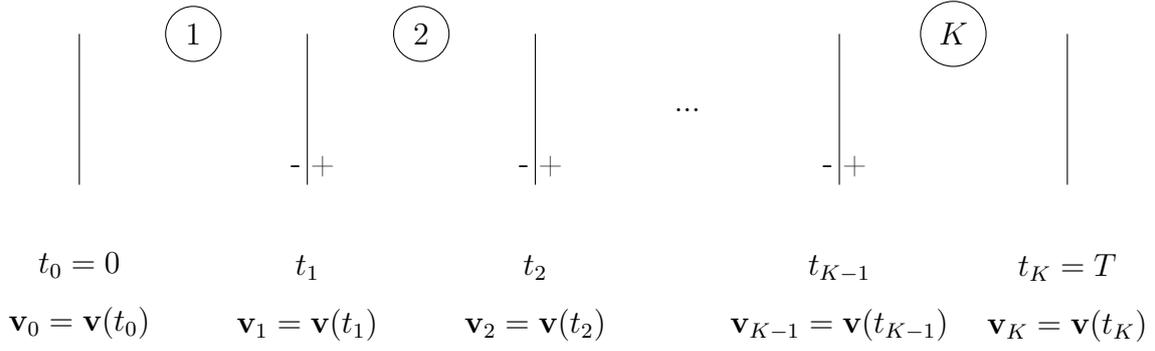
\begin{figure}[!htb]
\centering
\begin{tikzpicture}[node distance=2cm]
\draw (-9.5,0) -- (-9.5,2);
\draw (-6.5,0) -- (-6.5,2);
\draw (-3.5,0) -- (-3.5,2);
\draw (0.5,0) -- (0.5,2);
\draw (3.5,0) -- (3.5,2);
\node (t0) [tblock,align=center,xshift=-9.5cm,yshift=-1.47cm] {$t_0=0$ \\[0.25em] $\textbf{v}_0=\textbf{v}(t_0)$};
\node (t1) [tblock,align=center,xshift=-6.5cm,yshift=-1.5cm] {$t_1$ \\[0.25em] $\textbf{v}_1=\textbf{v}(t_1)$};
\node (t2) [tblock,align=center,xshift=-3.5cm,yshift=-1.5cm] {$t_2$ \\[0.25em] $\textbf{v}_2=\textbf{v}(t_2)$};
\node (tp) [tblock,align=center,xshift=0.5cm,yshift=-1.5cm] {$t_{K-1}$ \\[0.25em] $\textbf{v}_{K-1}=\textbf{v}(t_{K-1})$};
\node (tT) [tblock,align=center,xshift=3.5cm,yshift=-1.5cm] {$t_K=T$ \\[0.25em] $\textbf{v}_K=\textbf{v}(t_K)$};
\node (dots) [tblock, align=center, xshift=-1.5cm,yshift=1cm] {...};
\node (1) [tblock, align=center, xshift=0.70cm,yshift=0.25cm] {+};
\node (2) [tblock, align=center, xshift=0.35cm,yshift=0.25cm] {-};
\node (3) [tblock, align=center, xshift=-3.30cm,yshift=0.25cm] {+};
\node (4) [tblock, align=center, xshift=-3.65cm,yshift=0.25cm] {-};
\node (5) [tblock, align=center, xshift=-6.30cm,yshift=0.25cm] {+};
\node (6) [tblock, align=center, xshift=-6.65cm,yshift=0.25cm] {-};
\node at (-8,2) [circle,draw] (1) {$1$};
\node at (-5,2) [circle,draw] (2) {$2$};
\node at (2,2) [circle,draw] (P) {$K$};
\end{tikzpicture}
\caption{A sketch illustrating the segmenting approach of MSS. The constraint equations are propagated forward in time in all $K$ segments, such that $\textbf{v}(t_i^+)=\textbf{v}(t_i^-)$ is satisfied for $i=1,2,...,K-1$.}
\label{multiseg diag}
\end{figure}

For completeness, we provide below a short description of the solution of system \eqref{MSS eqn}; for all derivation details, the reader is referred to \cite{Blonigan2018}. The solution $\textbf{v}(t)$ that satisfies (\ref{MSS eqn}c,d) in all segments can be expressed as
\begin{equation}
\textbf{v}(t)=P_{t} \left(\phi^{t_i,t}\textbf{v}(t_i)+\int_{t_i}^t \phi^{\tau,t}\frac{\partial \textbf{f}}{\partial s} \, d\tau \right), \, \, \, \, \, \, \, \, \, \,  t_i \leq t < t_{i+1}
\end{equation}
where $\phi^{\tau,t}$ is the state transition matrix that satisfies 
\begin{equation}
\frac{d\phi^{t,\tau}}{d\tau}=\left(\frac{\partial f}{\partial u}\bigg|_\tau\right) \phi^{t,\tau}
\end{equation}
and $P_t$ is the projection operator (in discrete form it's a matrix) 
\begin{equation}
\label{projector eqn}
P_t=I-\frac{\textbf{f}(t)\textbf{f}(t)^T}{\textbf{f}(t)^T\textbf{f}(t)}
\end{equation}
that eliminates $\eta$ from (\ref{MSS eqn}c) and enforces (\ref{MSS eqn}d). This allows us to integrate (\ref{MSS eqn}c) without the time dilation term $\eta \textbf{f}$, i.e.\
\begin{equation}\label{cons proj eqn}
\frac{d\textbf{v}'}{dt}-\frac{\partial f}{\partial u}\textbf{v}'-\frac{\partial \textbf{f}}{\partial s} =0
\end{equation}
and recover $\textbf{v}(t)$ from $\textbf{v}(t)=P_t\textbf{v}'(t)$. 

The minimisation problem (\ref{MSS eqn}) can be written as a least squares problem:
\begin{subequations}\label{MSS system eqn}
\begin{align}
\underset{\textbf{v}_i}{\text{Minimise}}  &\, \, \, \frac{1}{2} \sum_{i=0}^K \|\textbf{v}_i\|_2^2 \\
\text{subject} \, \text{to} & \, \, \, \textbf{v}_{i+1} = \Phi_{i+1}\textbf{v}_i+\textbf{b}_{i+1}
\end{align}
\end{subequations}
where $\textbf{v}_i=\textbf{v}(t_i^+)$, $\Phi_{i+1}=P_{t_{i+1}}\phi^{t_i,t_{i+1}}$ and $\textbf{b}_{i+1}=P_{t_{i+1}}\int_{t_i}^{t_{i+1}} (\phi^{\tau,t})\sfrac{\partial \textbf{f}}{\partial s} \, d\tau$. Equation (\ref{MSS system eqn}b) satisfies  the original constraints (\ref{MSS eqn}b,c,d). System \eqref{MSS system eqn} can be written in matrix form as 
\begin{subequations}\label{MSS system eqn1}
\begin{align}
\underset{\textbf{v}_i}{\text{Minimise}}  &\, \, \, \frac{1}{2}\sum_{i=0}^K \|\textbf{v}_i\|_2^2 \\
\text{subject} \, \text{to} & \, \, \, A\underline{\textbf{v}}=\underline{\textbf{b}}
\end{align}
\end{subequations}
where 
\begin{equation}\label{A eqn}
A = \begin{bmatrix}
-\Phi_1 & I \\
& -\Phi_2 & I \\
& & \ddots &\ddots \\
& & & -\Phi_K & I
\end{bmatrix} 
\, \, \, \, \, \underline{\textbf{v}}= \begin{bmatrix}
\textbf{v}_0 \\ \textbf{v}_1 \\ \vdots \\\textbf{v}_K
\end{bmatrix} \, \, \, \, \, \underline{\textbf{b}}=\begin{bmatrix}
\textbf{b}_1 \\ \textbf{b}_2 \\ \vdots \\\textbf{b}_K
\end{bmatrix}
\end{equation}

\noindent $A$ is a $NK \times N(K+1)$ matrix and $\underline{\textbf{v}}$ and $\underline{\textbf{b}}$  are vectors of length $N(K+1)$ and $NK$, respectively. Equation (\ref{MSS system eqn1}b) represents an under-determined system; according to (\ref{MSS system eqn1}a) we seek the solution with the minimum Euclidean norm. This a well known minimisation problem in linear algebra \cite{Press2007}. A set of discrete adjoint variables $\underline{\textbf{w}}=\begin{bmatrix}
\textbf{w}_1 & \textbf{w}_2 & \dots & \textbf{w}_K
\end{bmatrix}^T$ is introduced, and an optimality system is derived:
\begin{equation}\label{KKT eqn}
\begin{bmatrix}
-I & A^T \\
A & 0
\end{bmatrix}\begin{bmatrix}
\underline{\textbf{v}} \\ \underline{\textbf{w}}
\end{bmatrix}=
\begin{bmatrix}
0 \\ \underline{\textbf{b}}
\end{bmatrix}
\end{equation}
It is more efficient to solve a linear system with the Schur complement of (\ref{KKT eqn}), 
\begin{equation}\label{schur eqn}
S\underline{\textbf{w}} 
=\begin{bmatrix}
\Phi_1\Phi^T_1+I & -\Phi^T_2 \\
-\Phi_2 & \Phi_2\Phi^T_2+I & -\Phi^T_3 \\
& \ddots & \ddots & \ddots \\
& &-\Phi_K & \Phi_K\Phi^T_K+I
\end{bmatrix} \begin{bmatrix}
\textbf{w}_1 \\ \textbf{w}_2 \\ \vdots \\ \textbf{w}_K\end{bmatrix}
=
\underline{\textbf{b}},
\end{equation}
\noindent where $S=AA^T$, because the number of unknowns is more than halved. The solution $\underline{\textbf{w}}$ is substituted into \eqref{KKT eqn} to obtain $\underline{\textbf{v}}$. The Schur complement matrix $S$ is block tri-diagonal, symmetric and positive definite, and has size $NK\times NK$. 

Equation \eqref{schur eqn} can be solved iteratively by supplying matrix-vector products $S\underline{\textbf{w}}^{(m)}$ at each iteration $m$ to a Krylov subspace solver, such as Conjugate Gradient. All products can be computed by calling a time-stepper for the forward (or adjoint) equations with the  appropriate initial (or terminal) conditions. For example, the product $\Phi_{i}$ with an arbitrary vector $\textbf{z}_{i-1}$ requires forward integration of the homogeneous form of \eqref{cons proj eqn}, i.e.
\begin{equation}\label{con hom eqn}
\frac{d\textbf{v}'}{dt}-\frac{\partial f}{\partial u}\textbf{v}' =0
\end{equation}
\noindent with the initial condition $\textbf{v}'(t_{i-1}^+)=\textbf{z}_{i-1}$ until $t=t_i$. The projection operation is then applied, yielding $\Phi_i\textbf{z}_{i-1}=P_{t_i}\textbf{v}'(t_i^-)$.
Similarly, the product of the transpose matrix $\Phi_i^T$ with an arbitrary $\textbf{z}_i$ requires integration of the homogeneous adjoint equation 
\begin{equation}\label{adj eqn}
\frac{d\textbf{w}}{dt}+\frac{\partial f}{\partial u}^T\textbf{w} =0
\end{equation}
\noindent backwards in time with the terminal condition $\textbf{w}(t_i^-)=P_{t_{i}}\textbf{z}_i$ until $t=t_{i-1}$. The product is then  given by $\Phi_i^T\textbf{z}_i=\textbf{w}(t_{i-1}^+)$. A detailed algorithm to solve (\ref{schur eqn}) is available in \cite{Blonigan2018}. An adjoint algorithm for multiple control parameters is also available. The sensitivity $\sfrac{d\bar{J}}{ds}$ is obtained by integrating
\begin{equation}\label{MSS sens eqn}
\frac{d\bar{J}}{ds}=\frac{1}{T}\sum_{i=0}^{K-1} \int_{t_i}^{t_{i+1}} \bigg \langle\frac{\partial J}{\partial u}\bigg |_t, \textbf{v}' \bigg \rangle \, dt + \frac{1}{T}\sum_{i=0}^{K-1} \frac{\langle \textbf{f}_{i+1}, \textbf{v}'(t_{i+1})\rangle}{\|\textbf{f}_{i+1}\|_2^2}(\overline{J}-J_{i+1}) + \frac{\partial\bar{J}}{\partial s}
\end{equation}
Convergence of the iterative method is however slow without preconditioning of the system. We propose an efficient preconditioner in Section (\ref{Schur-PC}) below.

%% file: Schur-PC.tex
The convergence rate of iterative Krylov subspace solvers for symmetric, positive definite matrices (like $S$) depends on the distribution of the matrix eigenvalues \cite{Saad_2003}. For such systems, the eigenvalues are all positive and real. If all of them are tightly clustered around a few single points away from the origin, then one would expect fast convergence. On the other hand, widely spread eigenvalues without tight clustering can lead to slow convergence. The objective of a preconditioner is to reduce the spread of the eigenvalues, and thereby reduce the condition number, $\kappa(S) = \frac{\mu_{max} (S)}{ \mu_{min} (S)}$, where $\mu_{max}(S)$ and $\mu_{min}(S)$ are the maximum and minimum eigenvalues of $S$, respectively. 

An extensive survey of preconditioners for saddle point problems, such as the MSS KKT system \eqref{KKT eqn}, is available in \cite{Benzi2005}. Preconditioners can be applied to the $2 \times 2$ block system \eqref{KKT eqn} or directly to the Schur complement system \eqref{schur eqn}. In either case, an easily invertible approximation of $S$ is required. There is an additional restriction, namely that the preconditioner should be matrix free, i.e.\ it should rely on matrix-vector products only (computing and storing $S$ is out of the question for long trajectories and large $N$).

Efficient approximations can be made if one takes into account the structure of the problem. For example, in the finite element solution of the incompressible Navier-Stokes equations, the Schur complement can be interpreted as a discretisation of a second order diffusion operator. This is not surprising because pressure, that plays the role of a Lagrange multiplier that enforces incompressibility, is governed by a Poisson equation. For this operator, the action of $S^{-1}$ can be efficiently approximated by a multigrid iteration (more details can be found in \cite{wathen_2015,ELMAN2014}).

For the present case, however, such a route is not viable. In \cite{Blonigan2014} a second order partial differential equation for the adjoint variable $\textbf{w}$ is derived, but this is too complex to solve, and it is not evident how it can be simplified in order to aid the construction of a preconditioner.

Recently, McDonald et al. \cite{McDonald2018PreconditioningEquations} proposed a block circulant preconditioner for the `all-at-once' evolution of a linear system of ODEs with constant coefficients. Here `all-at-once' means that the space/time problem is written as a monolithic linear system. The authors exploit the block Toeplitz structure of the system to develop an efficient preconditioner that results in the number of Krylov iterations being independent of the number of time-steps. Unfortunately, this approach is also unsuitable for our problem. The block Toeplitz structure of the matrix in \cite{McDonald2018PreconditioningEquations}  is a direct consequence of the fact that the coefficients of the ODE are constant. In our problem, however, the blocks $\Phi_i$ of matrix $A$ in equation \eqref{A eqn} depend on time $t$ and the position along the reference trajectory.   

In the course of this work, we tried several preconditioners that have been proposed in the literature for the solution of general saddle systems, which however do not exploit the properties of the underlying physical problem. For example, Cao et al. \cite{Cao2005} proposed a splitting of the block matrix (\ref{KKT eqn}) into two matrices, one of which is used as a left preconditioner. The splitting contains an adjustable parameter, and theoretical analysis shows that the largest eigenvalue of the preconditioned system is one (independent of the value of the adjustable parameter), provided that the preconditioner is applied exactly. In practice, however, this cannot be achieved, and the preconditioner must be applied approximately. Although the method relies on matrix-vector products only, it is very time-consuming. The preconditioner of Golub et al. \cite{Golub2005} guarantees that the eigenvalues remain bounded within two intervals $[-1,(\sfrac{1-\sqrt[]{5})}{2}]$ and $[1,(\sfrac{1+\sqrt[]{5})}{2}]$. This tight clustering of eigenvalues ensures fast convergence, but constructing the preconditioner is very expensive. Standard preconditioners for general matrices, such as incomplete LU decomposition \cite{Saad_2003}, are also not suitable, and they also require storage of the matrix.  



A new approach is therefore warranted. The central idea is to identify the fastest growing modes of matrix $A$ and annihilate them. In this way, it is expected that it would be easier to minimise the norm  (\ref{MSS system eqn1}a) while preserving the continuity of $\textbf{v}(t)$ across segments. This can be achieved using partial singular value decomposition.

\subsection {A preconditioner based on partial singular value decomposition}
Recall that the Schur complement system (\ref{schur eqn}) is $AA^T\underline{\textbf{w}} =\underline{\textbf{b}} $. The Singular  Value Decomposition (SVD) of $A$ reads
\begin{equation}
A=U\Sigma V^T
\end{equation}
where $U$ is a $NK\times NK$ unitary matrix, $\Sigma$ is a $NK\times N(K+1)$ quasi-diagonal matrix, and $V$ is a $N(K+1)\times N(K+1)$  also unitary matrix. The columns of $U$ contain the left singular vectors of $A$, while the right singular vectors of $A$ make up the columns of $V$. $\Sigma$ contains the singular values of $A$, which we denote by $\sigma(A)$, in the $NK\times NK$ diagonal sub-matrix (the last $N$ columns consist of zeros, which we ignore). The singular values are ordered from largest to smallest, i.e. the diagonal elements are $\Sigma_{ii}=\sigma_i \quad (i=1\dots NK)$ and $\sigma_1>\sigma_2>...>\sigma_{NK}$. Note that $\sigma(A)=\sqrt{\mu(S)}$, where ${\mu(S)}$ are the eigenvalues of $S$.  Using the fact that $V^TV=I$ (since the columns of $V$ are orthonormal) we can write $S$ as
\begin{equation}\label{SVD_of_Sigma}
S=AA^T=U\Sigma \Sigma^T U^T=U\Sigma^2U^T
\end{equation}
and since $U^TU=I$ we have
\begin{equation}
\label{SVD Inv eqn}
S^{-1}=U\Sigma^{-2}U^T
\end{equation}
Equation \eqref{SVD Inv eqn} forms an exact preconditioner, which of course is not practical to compute. However, we can form an approximate preconditioner using the leading $l$ singular modes only, i.e. perform a partial SVD. $S$ can be written as 
\begin{equation}
S=\begin{bmatrix}U_1 & U_2
\end{bmatrix}\begin{bmatrix}
\Sigma^{2}_1 & 0 \\ 0 & \Sigma^{2}_2\end{bmatrix}\begin{bmatrix}U^T_1 \\ U^T_2
\end{bmatrix}=U_1\Sigma_1^{2}U_1^T+U_2\Sigma_2^{2}U_2^T
\end{equation}
\noindent where $U$ is partitioned as $U=\begin{bmatrix} U_1 & U_2\end{bmatrix}$, with $U_1=\begin{bmatrix} u_1  & u_2 &...& u_l\end{bmatrix}$ 
and \\ $U_2=\begin{bmatrix} u_{l+1}  & u_{l+2} &...&  u_{NK}\end{bmatrix}$. Matrix $\Sigma_1^2=diag(\sigma_1^2,\sigma_2^2,...,\sigma_l^2)$ contains the singular values corresponding to $U_1$ and  $\Sigma_2^2=diag(\sigma_{l+1}^2,\sigma_{l+2}^2,...,\sigma_{NK}^2)$ contains the rest. Replacing the diagonal submatrix $\Sigma^{2}_2$ with the identity matrix $I_2=I_{(NK-l)}$, we get an approximation $\hat{S}_{(l)}$ 
\begin{equation}
\hat{S}_{(l)}=\begin{bmatrix}U_1 & U_2
\end{bmatrix}\begin{bmatrix}
\Sigma^{2}_1 & 0 \\ 0 & I_2\end{bmatrix}\begin{bmatrix}U^T_1 \\ U^T_2
\end{bmatrix}=U_1\Sigma_1^{2}U_1^T+U_2U_2^T
\end{equation}
which can be easily inverted to give
\begin{equation}\label{SVD Inv 2 eqn}
\hat{S}_{(l)}^{-1}=M_{(l)}=\begin{bmatrix}U_1 & U_2
\end{bmatrix}\begin{bmatrix}
\Sigma^{-2}_1 & 0 \\ 0 & I_2\end{bmatrix}\begin{bmatrix}U^T_1 \\ U^T_2
\end{bmatrix}=U_1\Sigma_1^{-2}U_1^T+U_2U_2^T
\end{equation}
The product $M_{(l)}S$ now becomes
\begin{equation}
M_{(l)}S= \begin{bmatrix}U_1 & U_2
\end{bmatrix}\begin{bmatrix}
I & 0 \\ 0 &\Sigma^2_2\end{bmatrix}\begin{bmatrix}U^T_1 \\ U^T_2
\end{bmatrix}
\end{equation}
\noindent which indicates that $M_{(l)}$ has deflated the $l$ largest  singular values of matrix $S$ to $1$, while leaving the rest unaltered. In order to avoid computing the columns $U_2$ in \eqref{SVD Inv 2 eqn}, we invoke the orthogonality relation $U_2U_2^T=I-U_1U_1^T$, and we get
\begin{equation}\label{ext M eqn}
M_{(l)}=U_1\Sigma_1^{-2}U_1^T+(I-U_1U_1^T)
\end{equation}
Therefore in order to form $M_{(l)}$, we need $\Sigma_1$ and $U_1$ only. We deploy $M_{(l)}$ as a left preconditioner to convert (\ref{schur eqn}) into a better conditioned system of the form 
\begin{equation}\label{PC1 sys eqn}
M_{(l)}S\underline{\textbf{w}} =M_{(l)}\underline{\textbf{b}}
\end{equation}
and we solve for $\underline{\textbf{w}}$. 

Preconditioning based on partial singular value decomposition has been applied in the past, for example in \cite{Rezghi2010} to solve ill-conditioned least squares problems arising in image de-blurring applications. A preconditioner based on the partial Krylov-Schur decomposition of a matrix (see \cite{Stewart_2002} for details) was also used in \cite{Waugh2013} to accelerate the computation of limit cycles for thermoacoustic systems. The preconditioner had a form very similar to  \eqref{ext M eqn}, the difference being that instead of the diagonal matrix $\Sigma_1^{-2}$ in the first term on the right-hand side, the inverse of an upper triangular $l\times l$ matrix was taken. The preconditioner reduced the condition number significantly, speeded up the convergence of GMRES, and lead to modest overall cost savings (mainly due to the cost of converging the eigenvalues used to form the preconditioner). S\'anchez and Net \cite{SANCHEZ2010} also used a similar preconditioner for accelerating the convergence of a multiple shooting algorithm for finding periodic orbits.

The Lanczos bidiagonalization algorithm \cite{Baglama2005} (implemented in the `svds' command of MATLAB) can be employed to form an approximation to $M_{(l)}$ using $q$ iterations. The algorithm  requires matrix-vector products only. The approximation of $M_{(l)}$ after $q$ iterations is
\begin{equation}
M_{(l)}^{(q)}=U_1^{(q)}\Sigma_1^{-2(q)}U_1^{T(q)}+(I-U^{(q)}_1U_1^{T(q)})
\end{equation}
Computing $U_1^{(q)}$ and $\Sigma_1^{(q)}$ using the function `svds' of MATLAB requires evaluation of a large number of matrix-vector products $Ax$ and $A^Tz$, where $x$ and $z$ are algorithm-generated vectors of length $N(K+1)$ and $NK$ respectively. More specifically, each iteration requires at least $l+2$ (the smallest permissible Krylov subspace dimension) applications of $A$ and $A^T$. Even with one or two iterations of the algorithm, the cost of computing $M_{(l)}^{(q)}$ would most likely outweigh the cost savings due to preconditioning. A cheaper alternative is therefore required. 

\subsection {A simplified block diagonal preconditioner}
We can form a much cheaper preconditioner if we consider matrix $\tilde{A}$, obtained from $A$ by neglecting the upper diagonal, i.e. 
\begin{equation} \label{approximate A}
\tilde{A} = \begin{bmatrix}
-\Phi_1 & 0 \\
& -\Phi_2 & 0 \\
& & \ddots &\ddots \\
& & & -\Phi_K & 0
\end{bmatrix} 
\end{equation}
The corresponding Schur complement  $\tilde{S}=\tilde{A}\tilde{A}^T$ is a block diagonal matrix that takes the form
\begin{equation}\label{approx schur eqn}
\tilde{S} 
=\begin{bmatrix}
\Phi_1\Phi^T_1 &   \\
  & \Phi_2\Phi^T_2 &  \\
&  & \ddots & \\
& &  & \Phi_K\Phi^T_K
\end{bmatrix}
\end{equation}
This form indicates that each segment is now decoupled from the two neighbouring segments on either side. A Block Diagonal Preconditioner (BDP), $\mathbf{M}_{BD}$, can be constructed that approximates the inverse, $\tilde{S}^{-1}$, i.e.
\begin{equation}
\label{def_MBD}
\mathbf{M_{BD}}\approx \tilde{S}^{-1}= \begin{bmatrix}
(\Phi_1\Phi^T_1)^{-1}  \\
& (\Phi_2\Phi^T_2)^{-1} &  \\
&  & \ddots &  \\
& & & (\Phi_K\Phi^T_K)^{-1}
\end{bmatrix} 
\end{equation}
Each diagonal block in \eqref{def_MBD} can be approximated as before using partial singular value decompositions i.e.
\begin{subequations}\label{general BDP eqn}
\begin{align}
&\mathbf{{M}_{BD}}^{(q)}_{(l)} = diag(\mathbf{M}^{(q)}_{(l),1},\mathbf{M}^{(q)}_{(l),2},...,\mathbf{M}^{(q)}_{(l),K}) \\
&\mathbf{M}^{(q)}_{(l),i}
={\underline{U}_{1,i}}{\underline{\Sigma}_{1,i}^{-2}}{\underline{U}_{1,i}^{T}}+(I-{\underline{U}_{1,i}}{\underline{U}_{1,i}^{T}}) 
\end{align}
\end{subequations}
\noindent where $i$ is the segment number, $q$ is the number of iterations and $l$ is the number of retained singular modes in each segment. ${\underline{\Sigma}_{1,i}}$ and $\underline{U}_{1,i}$ are the matrices corresponding to the $l$ singular values and the left singular vectors of $\Phi_i$, respectively. The superscript $(q)$ and subscript $(l)$ have been removed from ${\underline{\Sigma}_{1,i}}$ and $\underline{U}_{1,i}$ for clarity. 

The preconditioner $\mathbf{{M}_{BD}}^{(q)}_{(l)}$ has several advantages. Firstly, it is much cheaper to construct than $M^{(q)}_{(l)}$; in fact, it is $O(K)$ times cheaper. Secondly, its computation is fully parallelizable in time. Each processor is assigned to one segment and computations can proceed independently, because each segment is treated separately, i.e.\ there is no message passing between processors. Thirdly, only $l$ vectors $\underline{U}_{1,i}$ need storage in each segment. Finally, both $S$ and $\mathbf{{M}_{BD}}^{(q)}_{(l)}$ are symmetric, positive definite  matrices, making the Conjugate Gradient method applicable. 

We use $\mathbf{{M}_{BD}}^{(q)}_{(l)}$ as a left preconditioner for the original system, i.e. we solve
\begin{equation}\label{PC2 sys eqn}
\mathbf{{M}_{BD}}^{(q)}_{(l)}S\underline{\textbf{w}} =\mathbf{{M}_{BD}}^{(q)}_{(l)}\underline{\textbf{b}}
\end{equation}
In the following two sections, we investigate the condition number of the matrix $\mathbf{{M}_{BD}}^{(q)}_{(l)} S$ and the performance of the preconditioner for two standard problems; the Lorenz system and the Kuramoto Sivashinsky equation.

%% file: PC-Lorenz.tex
The well-known Lorenz system takes the form, 
\begin{equation}
\begin{aligned}
&\frac{dx}{dt}=\sigma(y-x) \\
&\frac{dy}{dt}=x(\rho - z) -y \\
&\frac{dz}{dt}=xy - \beta z
\end{aligned}
\end{equation}
where $\sigma$, $\rho$ and $\beta$ are system parameters. The Lorenz system is a common test case for chaotic sensitivity analysis applications. For $\sigma=10$ and $\beta =\sfrac{8}{3}$, the largest Lyapunov exponent $\lambda_{max}$ increases from $\lambda_{max}\approx 0.8$ at $\rho=24$ to $\lambda_{max}\approx1.7$ at $\rho=96$ \cite{Fryland1984Lyapunov-exponentModel} (roughly a linear growth with $\rho$ with frequent dips). The remaining exponents are $\lambda_2=0$ and $\lambda_3<0$. The objective function considered is
\begin{equation}
\overline{J} = \frac{1}{T} \int_0^T z \, dt
\end{equation}
Using $\sigma=10$ and $\beta=\sfrac{8}{3}$, sensitivities were sought with respect to $\rho$, i.e. $\sfrac{d\overline{J}}{d\rho}$. MATLAB's variable step Runge-Kutta solver (ode45) was used to compute the  trajectory $\textbf{u}_{ref}(t)$ and to perform all constraint and adjoint integrations. 

\begin{figure}[!htb]
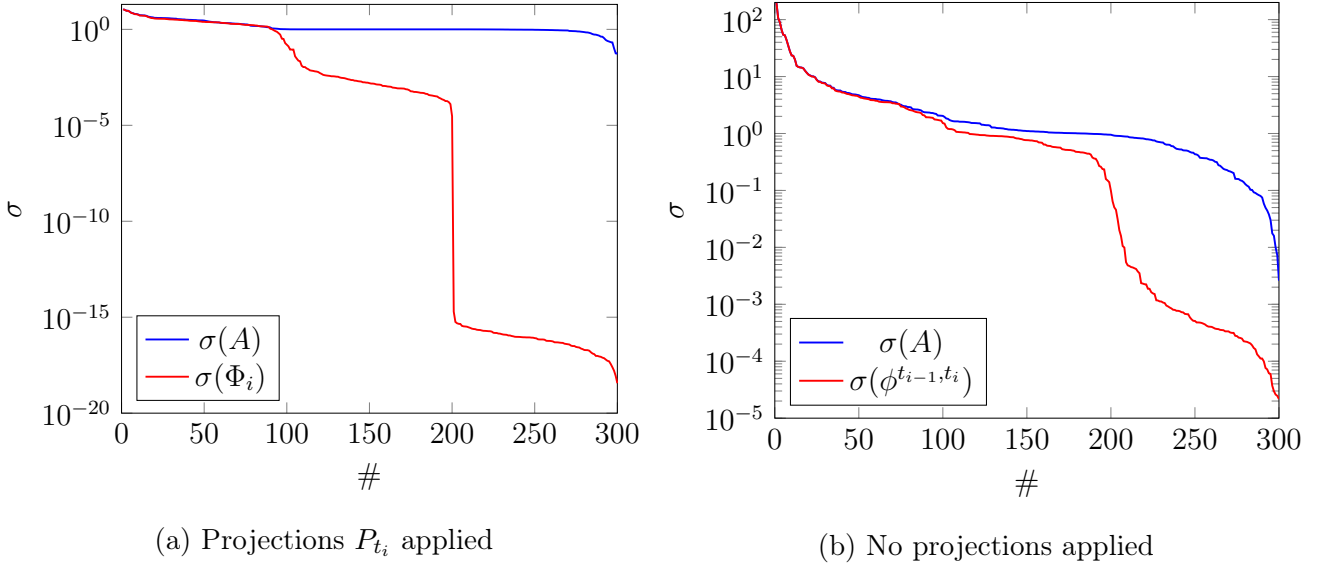

\begin{subfigure}{0.49\textwidth}
 \includestandalone[width=1\textwidth]{Sigma_Phi_A_general/Sigma_rho80_T50_dT05}
 \caption{Projections $P_{t_{i}}$ applied}
  \end{subfigure}
  \begin{subfigure}{0.49\textwidth}
\includestandalone[width=1\textwidth]{Sigma_Phi_A_general/Sigma_rho80_T50_dT05_np}
 \caption{No projections applied}
  \end{subfigure}
  \caption{A distribution of $\sigma(A)$, $\sigma(\Phi_i)=\sigma(P_{t_{i}}\phi^{t_{i-1},t_{i}})$ and $\sigma(\phi^{t_{i-1},t_{i}})$ ordered from the largest to the smallest values. Obtained for $\rho=80$ with $T=50$ and $\Delta T = 0.5$ ($K=100$ segments).}
\label{Sigma general fig}
\end{figure}

Figure (\ref{Sigma general fig}a) shows a comparison of the singular values of $A$ with the union of the singular values of all $\Phi_i$, ordered from largest to smallest. There are in total $3K=300$ singular values for the case examined. Each segment is located at a different place on the attractor, and has 3 local finite Lyapunov exponents associated with it. It is well known that such local exponents can fluctuate significantly around the average values as the trajectory is traced \cite{Wolfe_samelson_2007}. We observe that $\sigma(\Phi_i)$ is close to $\sigma(A)$ for the largest $K$ values, and thereafter the two curves start to deviate. The last $K$ values of $\sigma(\Phi_i)$ are very small,  which is expected, since $\Phi_{i}=P_{t_{i}}\phi^{t_{i-1},t_{i}}$ is almost singular due to the projection $P_{t_{i}}$. To remove the effect of singularity that distorts the comparison, Figure (\ref{Sigma general fig}b) shows $\sigma(A(\phi^{t_{i-1},t_{i}}))$ and $\sigma(\phi^{t_{i-1},t_{i}})$, where $A(\phi^{t_{i-1},t_{i}})$ is evaluated using $\phi^{t_{i-1},t_{i}}$ in equation \eqref{A eqn}, i.e.\ without applying the projection $P_{t_{i}}$ to the state transition matrix. The matching for the first $K$ singular values is now much more clearly seen.

In \cite{Blonigan2018} it was shown that $\mu_{max} (S)=\sigma_{max}^2(A)$ is related to the largest singular value of $\Phi_i$ across all segments, $\sigma_{max}(\Phi_i)$, as $\mu_{max}(S)=1+\sigma_{max}^2(\Phi_i)$, provided that $\sigma_{max}(\Phi_i) \gg 1$, in which case  $\mu_{max}(S) \approx \sigma_{max}^2(\Phi_i)$. Figure (\ref{Sigma general fig}) shows that this approximation holds for many more eigenvalues. For the Lorenz system, it holds for approximately K eigenvalues, each corresponding to the local positive Lyapunov exponent for each segment. 

\begin{figure}[!htb]
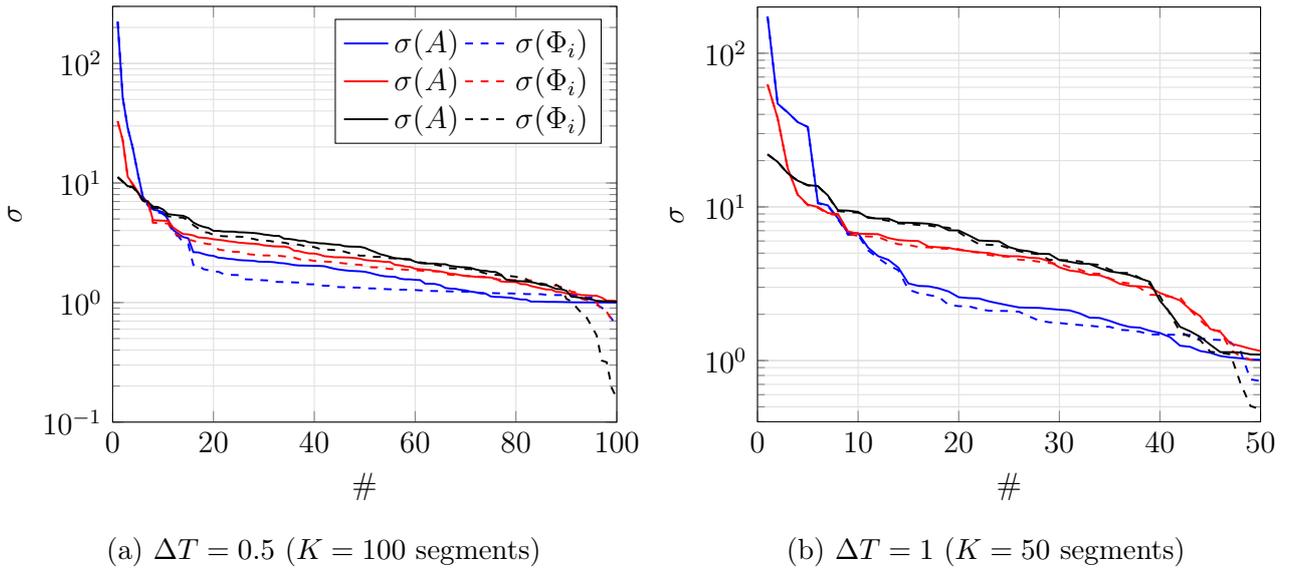

\begin{subfigure}{0.49\textwidth}
 \includestandalone[width=1\textwidth]{Sigma_Phi_A_dT05/sigma_phi_A_dT05}
 \caption{$\Delta T=0.5$ ($K=100$ segments)}
 \end{subfigure}
 \begin{subfigure}{0.49\textwidth}
 \includestandalone[width=0.96\textwidth]{Sigma_Phi_A_dT1/sigma_phi_A_dT1}
  \caption{$\Delta T=1$ ($K=50$ segments)}
 \end{subfigure}
  \caption{A distribution of the largest $K$ values of $\sigma(A)$ and the largest $\sigma(\Phi_i)$ of each segment. Obtained with $T=50$. Blue: $\rho=40$, red: $\rho=60$, black: $\rho=80$}
\label{Sigma dT fig}
\end{figure}

Figure (\ref{Sigma dT fig}) zooms in on the largest $K$ eigenvalues for two different $\Delta T$ values. The singular values  $\sigma(\Phi_i)$ approximate $\sigma(A)$ better for the larger $\Delta T$ value. This is because the diagonal blocks of matrix $A$ become more dominant for $\Delta T=1$ and the off-diagonal identity matrices can be neglected (refer to \eqref{approximate A}) without impairing the accuracy of the large singular values. 

\begin{figure}[t]
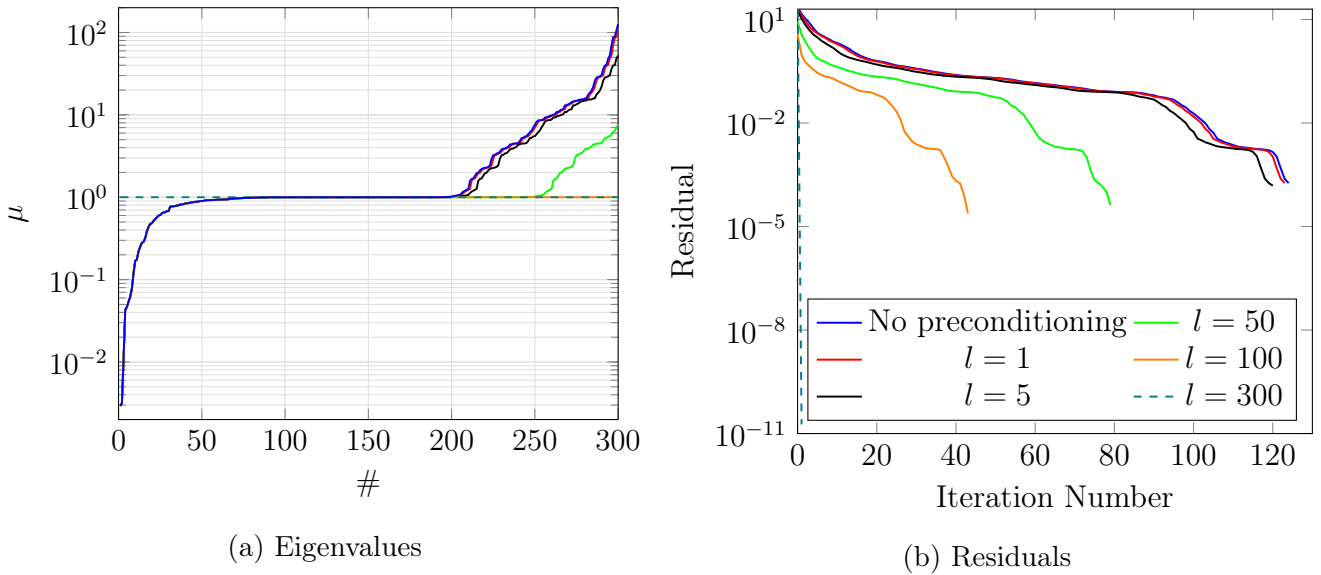

\begin{subfigure}{0.49\textwidth}
 \includestandalone[width=1\textwidth]{PC1/eigen}
         \caption{Eigenvalues}
\end{subfigure}
\begin{subfigure}{0.49\textwidth}
 \includestandalone[width=1\textwidth]{PC1/conv}
         \caption{Residuals}
\end{subfigure}
  \caption{Eigenvalues (ordered from smallest to largest) and convergence residuals for the original system $S$ (blue line) and the preconditioned system $M_{(l)}S$ for different $l$ (number of singular values). Obtained for $\rho=80$ with $T=50$ and $\Delta T = 0.5$ ($K=100$ segments).}
\label{PC1 fig}
\end{figure}

Next, we investigate the spectrum of the preconditioned system, starting with the exact preconditioner $M_{(l)}$ (\ref{ext M eqn}). The $l$ largest singular values and vectors of $S$ used to form $M_{(l)}$ have been obtained iteratively until convergence. Figure (\ref{PC1 fig}a) shows the eigenvalues of $\mu(S)$ and $\mu (M_{(l)}S)$. For the case examined, the matrix $S$ has approximately $K=100$ eigenvalues $\mu(S)>1$ (blue line), corresponding to the local positive Lyapunov exponent in each segment. Eigenvalues equal to $1$ correspond to the neutrally stable exponent, and the ones with $\mu(S)<1$ to the stable exponent. Different values of $l$ were used to construct $M_{(l)}$ (the values are reported in panel b). We notice that $M_{(l)}$ effectively deflates the $l$ largest eigenvalues of the preconditioned system $\mu(M_{(l)}S)$ to 1, while leaving the rest unaltered. As $l$ increases, more and more eigenvalues are clustered closer to 1, leading to increasingly faster convergence (shown in Figure \ref{PC1 fig}b). When $l=K(=100)$, $\mu_{max}(M_{(l)}S)=1$ and convergence still takes about $50$ iterations. For the limiting case of  $l=3K(=300)$ all eigenvalues are equal to $1$ (see left panel) and convergence is achieved in just one iteration, as expected. This indicates that we need not only deflate the large ($\mu$>1), but also to increase (regularise) the very small ($0<\mu \ll 1$) eigenvalues for fast convergence.  Such small eigenvalues make the system more singular, with implications to the accuracy of the computed sensitivity; this is explored later in Section (\ref{ill_cond}).  

\begin{figure}[htb!]
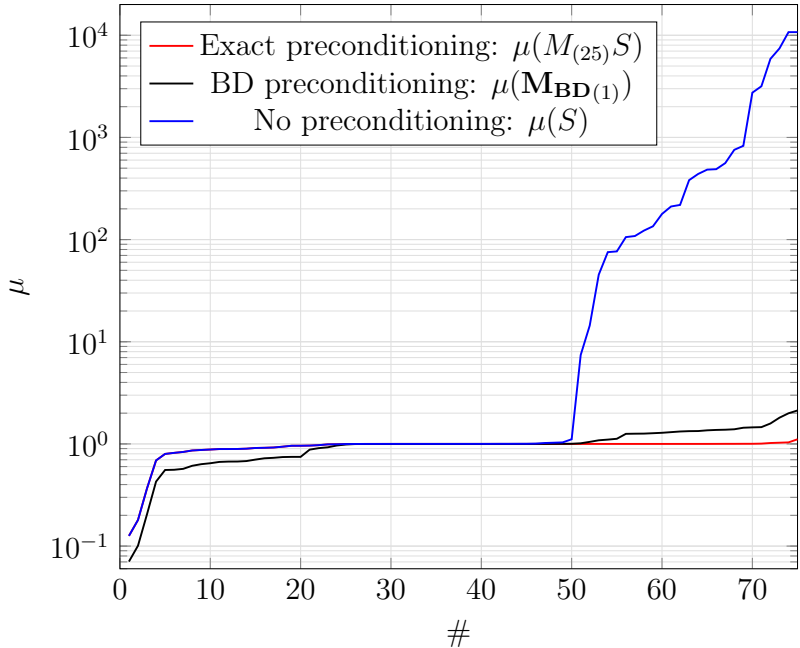

\centering
 \includestandalone[width=0.6\textwidth]{PC2_eigenvalues_T50/eigen_S_dT2}
  \caption{Eigenvalues of the original and preconditioned systems for $T=50$, $\Delta T=2$ ($K=25$) and $\rho=80$.}
\label{Eigen dT2 fig}
\end{figure}

We explore next the performance of the BDP \eqref{general BDP eqn}. Figure (\ref{Eigen dT2 fig}) shows $\mu(S)$, $\mu(M_{(25)}S)$ and $\mu(\mathbf{{M}_{BD}}_{(1)}S)$. The number of segments is $K=25$, and we have computed only one singular value until convergence in each segment. Matrix $\mathbf{{M}_{BD}}_{(1)}S$ has a very similar eigenvalue spectrum to $\mu(M_{(25)}S)$. It can be seen that $ \mu_{max}(M_{(25)}S)\approx 1$, while $ \mu_{max}(\mathbf{{M}_{BD}}_{(1)}S)\approx 2$. Most importantly, $\mathbf{{M}_{BD}}_{(1)}S$ has clustered the  $K$ largest eigenvalues in the interval $[1,2]$. This indicates that the approximate BDP  preconditioner encapsulates reliable information for the fastest growing modes, which results in the suppression of $ \mu_{max}(\mathbf{{M}_{BD}}_{(1)}S)$ by four orders of magnitude. However, there is a slight reduction in $\mu_{min}(\mathbf{{M}_{BD}}_{(1)}S)$ with respect to $\mu_{min}(S)$ and $\mu_{min}(M_{(25)}S)$. 

\begin{figure}[!htb]
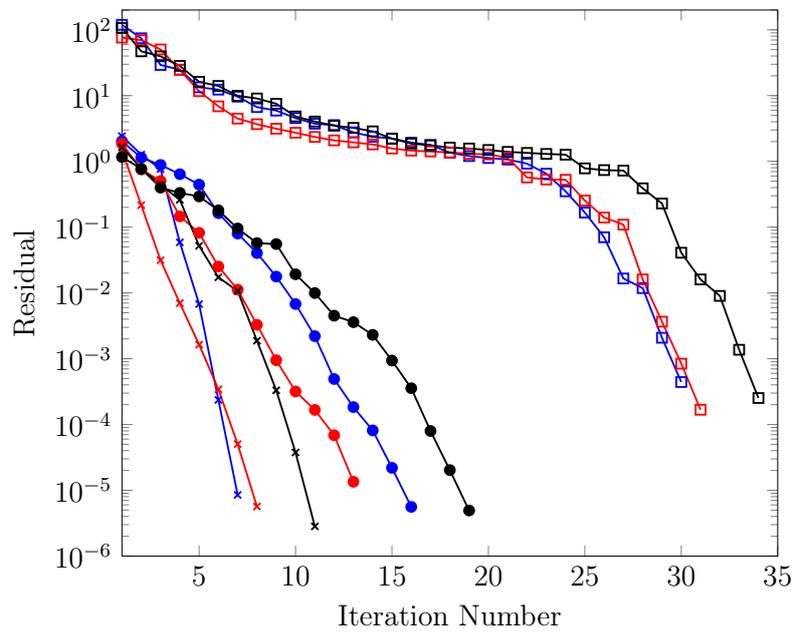

\centering
 \includestandalone[width=0.6\textwidth]{PC2_eigenvalues_T50/conv_S_dT2}
  \caption{Convergence history for the original and preconditioned systems using $T=50$ and $\Delta T=2$ ($K=25$). Blue: $\rho=40$, red: $\rho=60$, black: $\rho=80$. Squares: $S$, crosses: $M_{(25)}S$, circles: $\mathbf{{M}_{BD}}_{(1)}S$.}
\label{Conv dT2 fig}
\end{figure}

Figure (\ref{Conv dT2 fig}) shows the convergence rates for different $\rho$. As expected, using the exact preconditioner $M_{(l)}$ does provide considerably faster convergence compared to the original system (\ref{schur eqn}). Of course $M_{(25)}$ performs better than $\mathbf{{M}_{BD}}_{(1)}$ (since the latter assumes a block diagonal structure of $A$), but the cheaper cost of constructing and storing $\mathbf{{M}_{BD}}_{(1)}$ makes it a much more practical alternative.

%% file: Sigma_Phi_A_general/Sigma_rho80_T50_dT05.tex
\begin{tikzpicture}
\begin{semilogyaxis}[legend entries={$\sigma(A)$,$\sigma(\Phi_i)$},
  legend pos=south west,
	xlabel=$\#$,
  ylabel=$\sigma$,
   xmin=0,
   xmax=300,
   ymax=20,
   ymin=1e-20,
	no marks]
\addplot[color=blue,line width=0.75pt,solid,mark=x] %
table[x=No,y=sigma_A,col sep=comma]{Sigma_Phi_A_general/rho80.csv};
\addplot[color=red,line width=0.75pt,solid,mark=x] %
table[x=No,y=sigma_phi,col sep=comma]{Sigma_Phi_A_general/rho80.csv};
\end{semilogyaxis}
\end{tikzpicture}

%% file: Sigma_Phi_A_general/Sigma_rho80_T50_dT05_np.tex
\begin{tikzpicture}
\begin{semilogyaxis}[legend entries={$\sigma(A)$,$\sigma(\phi^{t_{i-1},t_{i}})$},
  legend pos=south west,
	xlabel=$\#$,
  ylabel=$\sigma$,
   xmin=0,
   xmax=300,
   ymax=200,
   ymin=1e-5,
	no marks]
\addplot[color=blue,line width=0.75pt,solid,mark=x] %
table[x=No,y=sigma_A_np,col sep=comma]{Sigma_Phi_A_general/rho80.csv};
\addplot[color=red,line width=0.75pt,solid,mark=x] %
table[x=No,y=sigma_phi_np,col sep=comma]{Sigma_Phi_A_general/rho80.csv};
\end{semilogyaxis}
\end{tikzpicture}

%% file: Sigma_Phi_A_dT05/sigma_phi_A_dT05.tex
\begin{tikzpicture}
\begin{semilogyaxis}
[legend entries={$\sigma(A)$,$\sigma(\Phi_i)$,$\sigma(A)$,$\sigma(\Phi_i)$,$\sigma(A)$,$\sigma(\Phi_i)$},
legend columns=2,
  legend pos=north east,
	xlabel=$\#$,
  ylabel=$\sigma$,
   xmin=0,
   xmax=100,
   ymax=300,
   ymin=1e-1,
	grid=both,
	minor grid style={gray!25},
	major grid style={gray!25},
	no marks]
\addplot[color=blue,line width=0.75pt,solid,mark=x]
table[x=No.,y=sigma_A,col sep=comma]{Sigma_Phi_A_dT05/rho40.csv};  \label{plot:line1}
\addplot[color=blue,line width=0.75pt,dashed,mark=x] 
table[x=No.,y=sigma_phi,col sep=comma]{Sigma_Phi_A_dT05/rho40.csv}; \label{plot:line2}
\addplot[color=red,line width=0.75pt,solid,mark=x] 
table[x=No.,y=sigma_A,col sep=comma]{Sigma_Phi_A_dT05/rho60.csv}; \label{plot:line3}
\addplot[color=red,line width=0.75pt,dashed,mark=x] 
table[x=No.,y=sigma_phi,col sep=comma]{Sigma_Phi_A_dT05/rho60.csv};
\addplot[color=black,line width=0.75pt,solid,mark=x] 
table[x=No.,y=sigma_A,col sep=comma]{Sigma_Phi_A_dT05/rho_80.csv};
\addplot[color=black,line width=0.75pt,dashed,mark=x] 
table[x=No.,y=sigma_phi,col sep=comma]{Sigma_Phi_A_dT05/rho_80.csv};
\end{semilogyaxis}
\end{tikzpicture}

%% file: Sigma_Phi_A_dT1/sigma_phi_A_dT1.tex
\begin{tikzpicture}
\begin{semilogyaxis}
[
  legend pos=south west,
	xlabel=$\#$,
  ylabel=$\sigma$,
   xmin=0,
   xmax=50,
   ymax=2e2,
   ymin=0.4,
	grid=both,
	minor grid style={gray!25},
	major grid style={gray!25},
	no marks]
\addplot[color=blue,line width=0.75pt,solid,mark=x]
table[x=No.,y=sigma_A,col sep=comma]{Sigma_Phi_A_dT1/rho40.csv};  
\addplot[color=blue,line width=0.75pt,dashed,mark=x] 
table[x=No.,y=sigma_phi,col sep=comma]{Sigma_Phi_A_dT1/rho40.csv};
\addplot[color=red,line width=0.75pt,solid,mark=x] 
table[x=No.,y=sigma_A,col sep=comma]{Sigma_Phi_A_dT1/rho60.csv}; 
\addplot[color=red,line width=0.75pt,dashed,mark=x] 
table[x=No.,y=sigma_phi,col sep=comma]{Sigma_Phi_A_dT1/rho60.csv};
\addplot[color=black,line width=0.75pt,solid,mark=x] 
table[x=No.,y=sigma_A,col sep=comma]{Sigma_Phi_A_dT1/rho80.csv};
\addplot[color=black,line width=0.75pt,dashed,mark=x] 
table[x=No.,y=sigma_phi,col sep=comma]{Sigma_Phi_A_dT1/rho80.csv};
 \coordinate (legend) at (axis description cs:0.97,0.03);
\end{semilogyaxis}
\end{tikzpicture}

%% file: PC1/eigen.tex
\begin{tikzpicture}
\begin{semilogyaxis}[
	xlabel=$\#$,
  ylabel=$\mu$,
   xmin=0,
   xmax=300,
   ymin=2e-3,
   ymax=200,
	grid=both,
	minor grid style={gray!25},
	major grid style={gray!25},
	no marks]
\addplot[color=red,line width=0.75pt,solid] 
table[x=No.,y=1,col sep=comma]{PC1/eigen.csv};  
\addplot[color=black,line width=0.75pt,solid] 
table[x=No.,y=5,col sep=comma]{PC1/eigen.csv};  
\addplot[color=green,line width=0.75pt,solid] 
table[x=No.,y=50,col sep=comma]{PC1/eigen.csv};  
\addplot[color=orange,line width=0.75pt,solid] 
table[x=No.,y=100,col sep=comma]{PC1/eigen.csv};  
\addplot[color=teal,line width=0.75pt,dashed] 
table[x=No.,y=300,col sep=comma]{PC1/eigen.csv}; 
\addplot[color=blue,line width=0.75pt,solid]
table[x=No.,y=No_PC,col sep=comma]{PC1/eigen.csv}; 
\end{semilogyaxis}
\end{tikzpicture}

%% file: PC1/conv.tex
\begin{tikzpicture}
\begin{semilogyaxis}
[
legend entries={No preconditioning,$l=50$,$l=1$,$l=100$,$l=5$,$l=300$},
legend columns=2,
  legend pos=south east,
	xlabel=Iteration Number,
  ylabel=Residual,
   xmin=0,
   xmax=130,
   ymax=20,
   ymin=1e-11,
	no marks]
\addplot[color=blue,line width=0.75pt,solid,mark=x]
table[x=x,y=No_PC,col sep=comma]{PC1/res.csv};  
\addplot[color=green,line width=0.75pt,solid,mark=x] 
table[x=x,y=50,col sep=comma]{PC1/res.csv};  
\addplot[color=red,line width=0.75pt,solid,mark=x] 
table[x=x,y=1,col sep=comma]{PC1/res.csv};  
\addplot[color=orange,line width=0.75pt,solid,mark=x] 
table[x=x,y=100,col sep=comma]{PC1/res.csv};  
\addplot[color=black,line width=0.75pt,solid,mark=x] 
table[x=x,y=5,col sep=comma]{PC1/res.csv};  
\addplot[color=teal,line width=0.75pt,dashed,mark=x] 
table[x=x,y=300,col sep=comma]{PC1/res.csv};  
\end{semilogyaxis}
\end{tikzpicture}

%% file: PC2_eigenvalues_T50/eigen_S_dT2.tex
\begin{tikzpicture}
\begin{semilogyaxis}
[legend entries={Exact preconditioning: $\mu(M_{(25)}S)$,BD preconditioning:  $\mu(\mathbf{{M}_{BD}}_{(1)})$,No preconditioning: $\mu(S)$},
  legend pos=north west,
	xlabel=$\#$,
  ylabel=$\mu$,
   xmin=0,
   xmax=75,
   ymax=2e4,
   ymin=6e-2,
	grid=both,
	minor grid style={gray!25},
	major grid style={gray!25},
	width=0.62\textwidth,
	no marks]
\addplot[color=red,line width=0.75pt,solid,mark=x] 
table[x=No.,y=eig_pc1,col sep=comma]{PC2_eigenvalues_T50/rho80.csv};
\addplot[color=black,line width=0.75pt,solid,mark=x] 
table[x=No.,y=eig_pc2,col sep=comma]{PC2_eigenvalues_T50/rho80.csv};
\addplot[color=blue,line width=0.75pt,solid,mark=x]
table[x=No.,y=eig_nopc,col sep=comma]{PC2_eigenvalues_T50/rho80.csv};   
\end{semilogyaxis}
\end{tikzpicture}

%% file: PC2_eigenvalues_T50/conv_S_dT2.tex
\begin{tikzpicture}
\begin{semilogyaxis}
[
	xlabel=Iteration Number,
  ylabel=Residual,
   xmin=1,
   xmax=35,
   ymax=2e2,
   ymin=1e-6,
	width=0.62\textwidth]
\addplot[color=blue,line width=0.75pt,solid,mark=square]
table[x=No.,y=conv_nopc,col sep=comma]{PC2_eigenvalues_T50/rho40.csv};  
\addplot[color=blue,line width=0.75pt,solid,mark=x] 
table[x=No.,y=conv_pc1,col sep=comma]{PC2_eigenvalues_T50/rho40.csv};
\addplot[color=blue,line width=0.75pt,solid,mark=*] 
table[x=No.,y=conv_pc2,col sep=comma]{PC2_eigenvalues_T50/rho40.csv}; 
\addplot[color=red,line width=0.75pt,solid,mark=square]
table[x=No.,y=conv_nopc,col sep=comma]{PC2_eigenvalues_T50/rho60.csv};  
\addplot[color=red,line width=0.75pt,solid,mark=x] 
table[x=No.,y=conv_pc1,col sep=comma]{PC2_eigenvalues_T50/rho60.csv};
\addplot[color=red,line width=0.75pt,solid,mark=*] 
table[x=No.,y=conv_pc2,col sep=comma]{PC2_eigenvalues_T50/rho60.csv}; 
\addplot[color=black,line width=0.75pt,solid,mark=square]
table[x=No.,y=conv_nopc,col sep=comma]{PC2_eigenvalues_T50/rho80.csv};  
\addplot[color=black,line width=0.75pt,solid,mark=x] 
table[x=No.,y=conv_pc1,col sep=comma]{PC2_eigenvalues_T50/rho80.csv};
\addplot[color=black,line width=0.75pt,solid,mark=*] 
table[x=No.,y=conv_pc2,col sep=comma]{PC2_eigenvalues_T50/rho80.csv}; 
\end{semilogyaxis}
\end{tikzpicture}

%% file: BDP_KS_App.tex
\begin{figure}[!htb]
\centering
\includegraphics{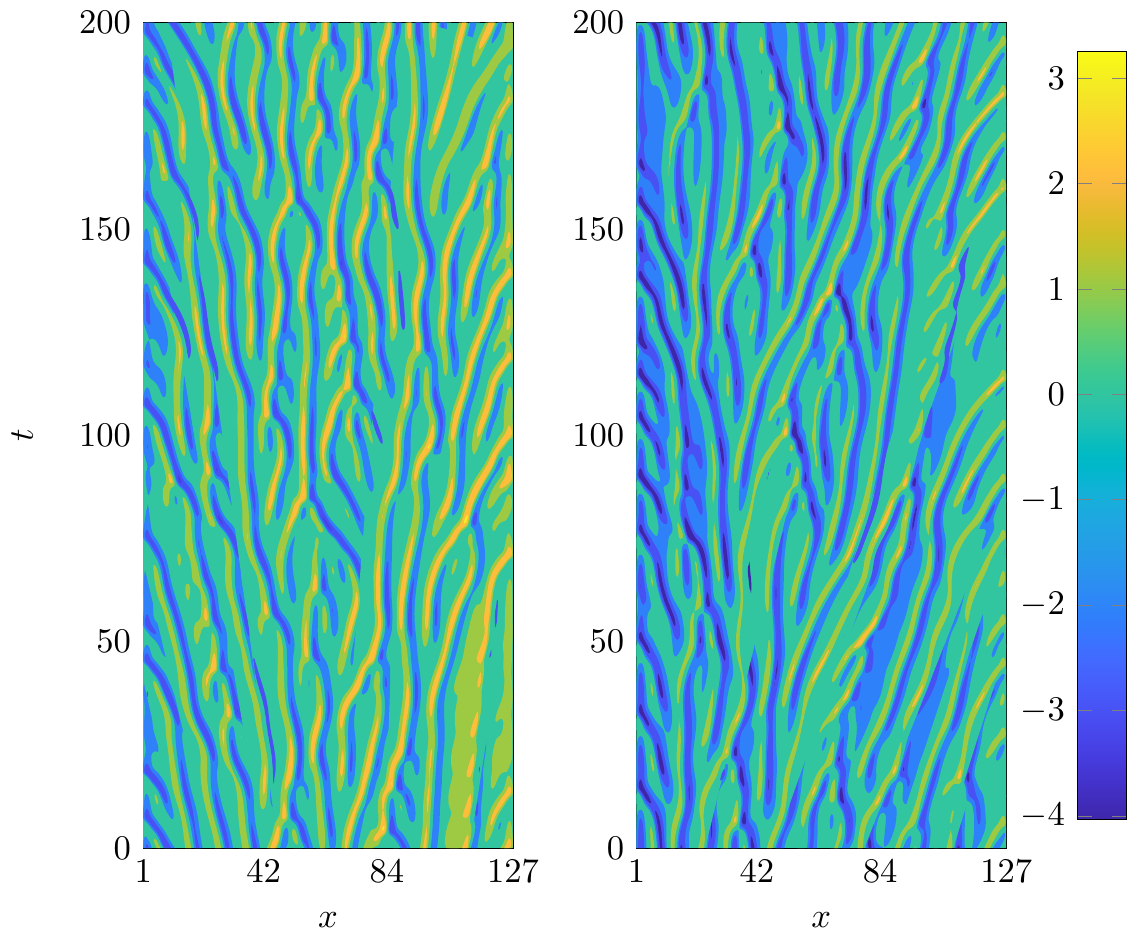}%
  \caption{Space-time plot of the solution $\textbf{u}(x,t)$ for $L=128$ using $N=255$ nodes in the $x$-direction (left: $c=0$, right: $c=0.8$). The integration time interval is $[-1000,200]$ with initial condition $\textbf{u}_0=1$.}
\label{KS solution fig}
\end{figure}

In this section we apply the MSS method with the BDP (\ref{general BDP eqn}) to a slightly modified version of the Kuramoto Sivashinsky (KS) equation \cite{Blonigan2014a}:
\begin{equation}\label{KS eqn}
\begin{aligned}
&\frac{\partial u}{\partial t}=-(u+c)\frac{\partial u}{\partial x}-\frac{\partial^2 u}{\partial x^2}-\frac{\partial^4 u}{\partial x^4} \\
&x\in[0,L] \\
& u(0,t)=u(L,t)=0 \\
&\frac{\partial u}{\partial x}\bigg |_{x=0}=\frac{\partial u}{\partial x}\bigg |_{x=L}=0
\end{aligned}
\end{equation}
where $L=128$ to ensure chaotic solutions \cite{Hyman1986TheSystems}. The Dirichlet and Neumann boundary conditions ensure ergodicity. The spatial derivatives were discretised into $N+2$ nodes ($N$ interior nodes and two boundary nodes) using second order central finite difference approximations on a uniform grid. The Neumann boundary conditions were enforced using ghost nodes. The MATLAB command ode45 was used for time integration. Figure (\ref{KS solution fig}) shows the solution to (\ref{KS eqn}) for $c=0$ and $c=0.8$. In both cases organised structures with a dominant wavelength appear \cite{holmes_lumley_berkooz_rowley_2012}, corresponding to the light turbulence regime.

\noindent Two objective functions were considered:
\begin{subequations}\label{KS cost functions eqn}
\begin{align}
\langle \bar{u} \rangle &=\frac{1}{TL}\int_0^T \int_0^L u  \, dx \, dt\\
\langle \overline{u^2} \rangle &=\frac{1}{TL}\int_0^T \int_0^L u^2  \, dx \, dt
\end{align}
\end{subequations}
and their sensitivities to parameter $c$, $\sfrac{d\langle \bar{u} \rangle}{dc}$ and $\sfrac{d \langle \overline{u^2} \rangle}{dc}$, were sought. The MSS segment size for all the cases studied was $\Delta T=10$ (based on $\lambda_{max}\approx 0.1$ for $c=0$ and $c=0.8$ \cite{Blonigan2014a}) unless otherwise stated.

The ideal use of the BDP (\ref{general BDP eqn}) involves choosing the parameters $q$ and $l$ such that the total number of matrix-vector products involving $\Phi_i$ and $\Phi^T_i$ is minimised. Considering that there are 15 positive Lyapunov exponents for $c=0.8$ \cite{Blonigan2014a}, it seems reasonable to choose $l=15$ in each segment to construct (\ref{general BDP eqn}). We explore the effect of different values of $l$ later.

\begin{figure}[!htb]
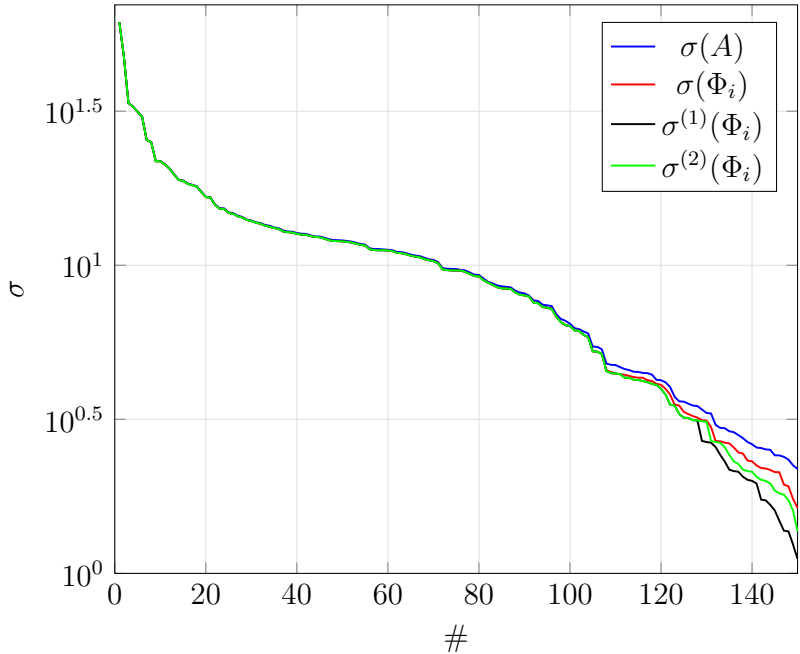

\centering
 \includestandalone[width=0.6\textwidth]{KS_effect_no_it/sigma}
  \caption{$\sigma(A)$ (blue line) and the largest 15 $\sigma(\Phi_i)$ for all segments ordered from largest to smallest (red: exact, black: $q=1$ iteration, green: $q=2$ iterations). Computed for $N=127$, $c=0.8$, $T=100$ ($K=10$).}
\label{sigma no it fig}
\end{figure}

Figure (\ref{sigma no it fig}) shows $\sigma(A)$ and the union of the computed 15 singular values of $\sigma(\Phi_i)$ for each segment, ordered from largest to smallest. Both the converged values and the values with $q=1$ and $q=2$ iterations are shown. It is clear that we can obtain accurate approximations to $\sigma(A)$ even with just 1 iteration when we evaluate $\sigma(\Phi_i)$. The curves start to deviate when the singular values approach unity. 

\begin{figure}[!htb]
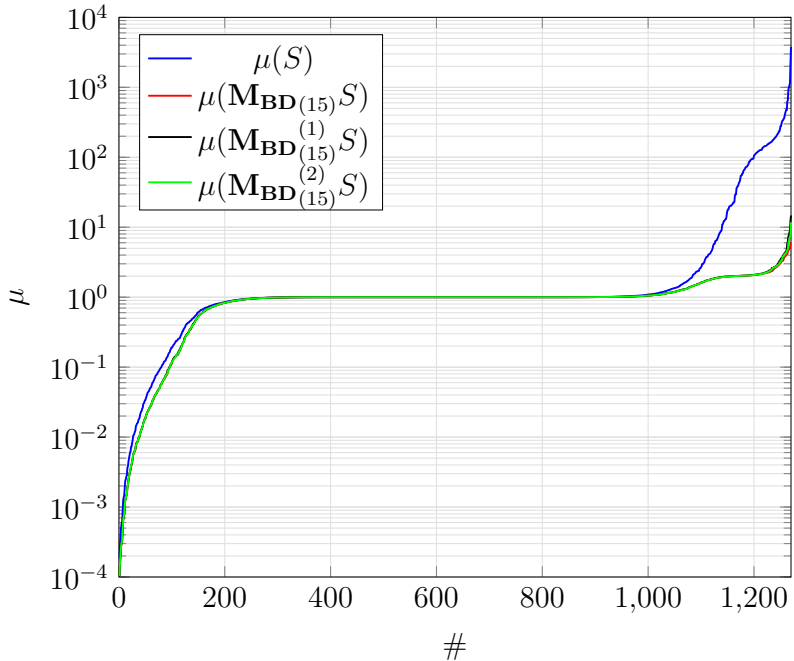

\centering
 \includestandalone[width=0.6\textwidth]{KS_effect_no_it/eig}
  \caption{Eigenvalues of the original system (blue line), the exact BDP ($\mathbf{{M}_{BD}}_{(15)}$) and inexact BDP ($\mathbf{{M}_{BD}}^{(1)}_{(15)}$,$\mathbf{{M}_{BD}}^{(2)}_{(15)}$) using $l=15$.}
\label{eig no it fig}
\end{figure}

Figure (\ref{eig no it fig}) shows the eigenvalues of the preconditioned system. Using $\mathbf{{M}_{BD}}^{(q)}_{(l)}$ has reduced the maximum eigenvalue (and therefore the condition number $\kappa$) by more than 2 orders of magnitude. The spectra of the inexact BDP systems  $\mathbf{{M}_{BD}}^{(1)}_{(15)}S$ and $\mathbf{{M}_{BD}}^{(2)}_{(15)}S$ are very similar to the spectra of the exact BDP system $\mathbf{{M}_{BD}}_{(15)}S$,  (with $\mu_{max}(\mathbf{{M}_{BD}}^{(1)}_{(15)})\approx 15$,  $\mu_{max}(\mathbf{{M}_{BD}}^{(2)}_{(15)}S)\approx 12$ and $\mu_{max}(\mathbf{{M}_{BD}}_{(15)}S)\approx 6$).
\begin{figure}[!b]
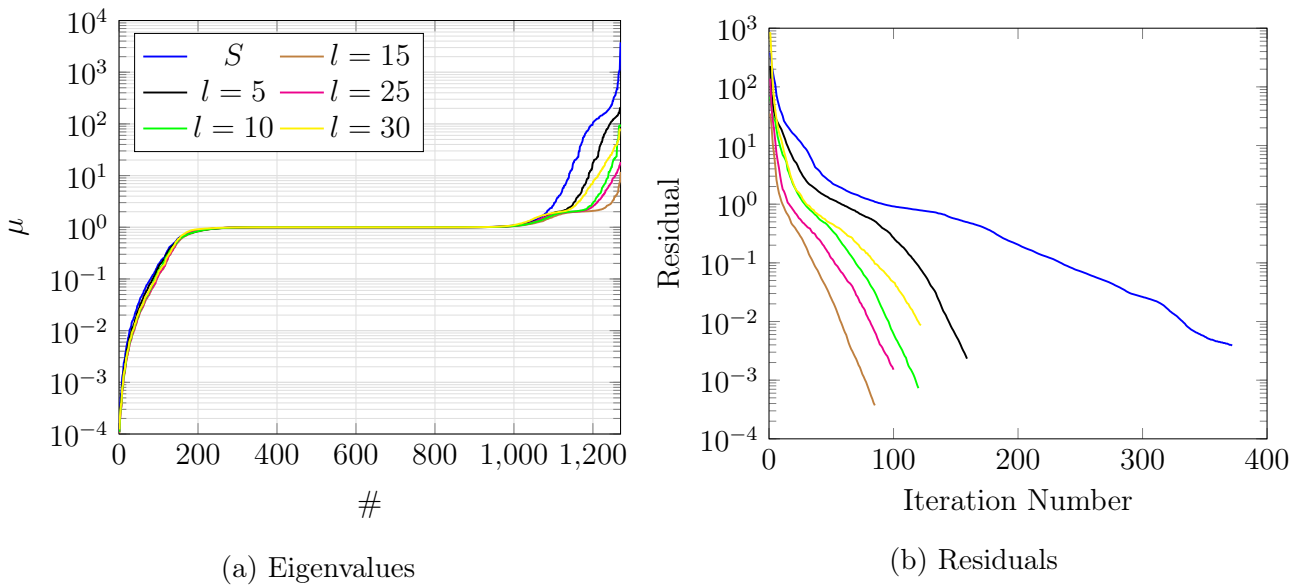

\begin{subfigure}{0.48\textwidth}
 \includestandalone[width=1\textwidth]{KS_effect_no_sigma/Eig}
\caption{Eigenvalues}
\end{subfigure}
\begin{subfigure}{0.49\textwidth}
 \includestandalone[width=1\textwidth]{KS_effect_no_sigma/res}
\caption{Residuals}
\end{subfigure}
\caption{Eigenvalues and residuals of the original system $S$ (blue line), and of the BDP system $(\mathbf{{M}_{BD}}^{(2)}_{(l)}S)$ for $N=127$, $T=100$ and $c=0.8$. The preconditioners were constructed for different $l$, and their residuals were found with a regularisation value $\gamma=0.01$. }
\label{eig no sigma fig}
\end{figure}

We study next the effect of $l$ on the spectra of  $\mathbf{{M}_{BD}}^{(q)}_{(l)}S$ and on the convergence rate.  The value $q=2$ was kept constant. Results for different $l$ are shown in Figure (\ref{eig no sigma fig}). Increasing $l$ up to $l=15$ improves the clustering of eigenvalues and leads to faster convergence rates (panel b). Interestingly, a further increase to $l=25$ or $l=30$ increases the condition number and slows down the convergence (Figure \ref{eig no sigma fig}b). This indicates that after a certain value of $l$, adding more singular modes starts to provide unreliable information to the preconditioner $\mathbf{{M}_{BD}}^{(q)}_{(l)}$. This is because as $l$ increases, the singular values approach unity, and  $\sigma(\Phi_i)$ starts to diverge from $\sigma(A)$. This was also observed clearly in Figure (\ref{Sigma general fig}) for the Lorenz system. 

\begin{figure}[!htb]
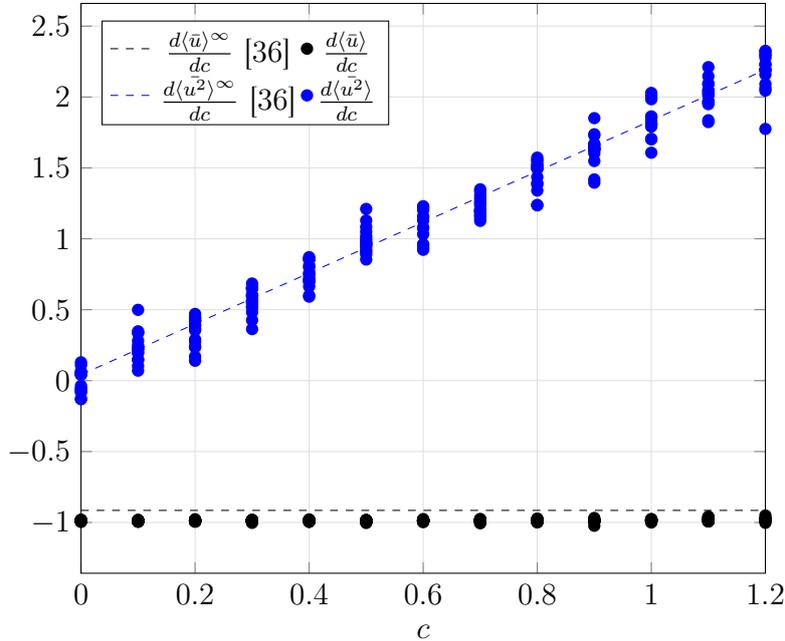

\centering
 \includestandalone[width=0.6\textwidth]{KS_validation/sensitivities}
  \caption{Sensitivities of $\langle \bar{u}\rangle$ and $\langle \overline{u^2}\rangle$ to parameter $c$. The dashed lines \cite{Blonigan2014a} are shown for reference and were obtained by differentiating curve fits for $T\to \infty$. The black dots ($\sfrac{d\langle \bar{u}\rangle}{dc}$) and the blue dots ($\sfrac{d\langle \overline{u^2}\rangle}{dc}$) were obtained for $T=100$ trajectories with random $\textbf{u}_0$ using MSS for $N=127$ and $N=255$, respectively.}
\label{KS sensitivities fig}
\end{figure}

Figure (\ref{KS sensitivities fig}) shows $\sfrac{d \langle \bar{u} \rangle}{dc}$ and $\sfrac{d \langle \overline{u^2} \rangle}{dc}$ for 20 and 15 different initial conditions $\textbf{u}_0$ respectively. The values of parameter $c$ examined are between $0 $ and $1.2$, which correspond to the light turbulence regime. Each data point was computed for a randomly generated initial condition vector $0<\textbf{u}_0<1$. To obtain $\textbf{u}(x,t)$, (\ref{KS eqn}) was integrated in the interval $[-1000,100]$, and MSS was applied to $\textbf{u}(x,t)$ in the interval $[0,100]$. The reference data \cite{Blonigan2014a} is the derivative of the curve fit of $\langle\bar{u}\rangle$ and $\langle \overline{u^2} \rangle$ vs. $c$, obtained for very long trajectories ($T=2000$). Figure (\ref{KS sensitivities fig}) shows that MSS slightly under-predicts $\sfrac{d\langle\bar{u}\rangle}{dc}$  (similar to \cite{Blonigan2014a,Blonigan2018}). This difference will be discussed in the next section.

%% file: KS_effect_no_it/sigma.tex
\begin{tikzpicture}
\begin{semilogyaxis}
[legend entries={$\sigma(A)$,$\sigma(\Phi_i)$,$\sigma^{(1)}(\Phi_i)$,$\sigma^{(2)}(\Phi_i)$},
  legend pos=north east,
	xlabel=$\#$,
  ylabel=$\sigma$,
   xmin=0,
   xmax=150,
   ymax=70,
   ymin=1,
	grid=both,
	minor grid style={gray!25},
	major grid style={gray!25},
	width=0.62\textwidth,
	no marks]
\addplot[color=blue,line width=0.75pt,solid,mark=x]
table[x=No.,y=Sigma(A)ex,col sep=comma]{KS_effect_no_it/sigma.csv};  
\addplot[color=red,line width=0.75pt,solid,mark=x] 
table[x=No.,y=Sigma(phi)ex,col sep=comma]{KS_effect_no_it/sigma.csv};  
\addplot[color=black,line width=0.75pt,solid,mark=x] 
table[x=No.,y=Sigma(phi)1,col sep=comma]{KS_effect_no_it/sigma.csv};  
\addplot[color=green,line width=0.75pt,solid,mark=x] 
table[x=No.,y=Sigma(phi)2,col sep=comma]{KS_effect_no_it/sigma.csv};  
\end{semilogyaxis}
\end{tikzpicture}

%% file: KS_effect_no_it/eig.tex
\begin{tikzpicture}
\begin{semilogyaxis}
[legend entries={$\mu(S)$,$\mu(\mathbf{{M}_{BD}}_{(15)}S)$,$\mu(\mathbf{{M}_{BD}}^{(1)}_{(15)}S)$,$\mu(\mathbf{{M}_{BD}}^{(2)}_{(15)}S)$},
  legend pos=north west,
	xlabel=$\#$,
  ylabel=$\mu$,
   xmin=0,
   xmax=1270,
   ymax=1e4,
   ymin=1e-4,
	grid=both,
	minor grid style={gray!25},
	major grid style={gray!25},
	width=0.62\textwidth,
	no marks]
\addplot[color=blue,line width=0.75pt,solid,mark=x]
table[x=no.eig,y=eig_no_pc,col sep=comma]{KS_effect_no_it/eig_res.csv};  
\addplot[color=red,line width=0.75pt,solid,mark=x] 
table[x=no.eig,y=eig_ex,col sep=comma]{KS_effect_no_it/eig_res.csv};  
\addplot[color=black,line width=0.75pt,solid,mark=x] 
table[x=no.eig,y=eig_1,col sep=comma]{KS_effect_no_it/eig_res.csv};  
\addplot[color=green,line width=0.75pt,solid,mark=x] 
table[x=no.eig,y=eig_2,col sep=comma]{KS_effect_no_it/eig_res.csv};  
\end{semilogyaxis}
\end{tikzpicture}

%% file: KS_effect_no_sigma/Eig.tex
\begin{tikzpicture}
\begin{semilogyaxis}
[legend entries={$S$,$l=15$,$l=5$,$l=25$,$l=10$,$l=30$},
legend columns=2,
ytick={1e-4,1e-3, 1e-2, 1e-1,1e0,1e1,1e2,1e3,1e4},
  legend pos=north west,
	xlabel=$\#$,
   ylabel=$\mu$,
   xmin=0,
   xmax=1270,
   ymax=1e4,
   ymin=1e-4,
	grid=both,
	minor grid style={gray!25},
	major grid style={gray!25},
	no marks]
\addplot[color=blue,line width=0.75pt,solid,mark=x]
table[x=no.,y=no_pc,col sep=comma]{KS_effect_no_sigma/eig.csv};  
\addplot[color=brown,line width=0.75pt,solid,mark=x] 
table[x=no.,y=15,col sep=comma]{KS_effect_no_sigma/eig.csv};  
\addplot[color=black,line width=0.75pt,solid,mark=x] 
table[x=no.,y=5,col sep=comma]{KS_effect_no_sigma/eig.csv};  
\addplot[color=magenta,line width=0.75pt,solid,mark=x] 
table[x=no.,y=25,col sep=comma]{KS_effect_no_sigma/eig.csv};  
\addplot[color=green,line width=0.75pt,solid,mark=x] 
table[x=no.,y=10,col sep=comma]{KS_effect_no_sigma/eig.csv};  
\addplot[color=yellow,line width=0.75pt,solid,mark=x] 
table[x=no.,y=30,col sep=comma]{KS_effect_no_sigma/eig.csv};  
\end{semilogyaxis}
\end{tikzpicture}

%% file: KS_effect_no_sigma/res.tex
\begin{tikzpicture}
\begin{semilogyaxis}
[
  legend pos=north west,
	xlabel=Iteration Number,
   ylabel=Residual,
   xmin=0,
   xmax=400,
   ymax=1e3,
   ymin=1e-4,
	no marks]
\addplot[color=blue,line width=0.75pt,solid,mark=x]
table[x=no.,y=no_pc,col sep=comma]{KS_effect_no_sigma/res.csv};  
\addplot[color=black,line width=0.75pt,solid,mark=x] 
table[x=no.,y=5,col sep=comma]{KS_effect_no_sigma/res.csv};  
\addplot[color=green,line width=0.75pt,solid,mark=x] 
table[x=no.,y=10,col sep=comma]{KS_effect_no_sigma/res.csv};  
\addplot[color=brown,line width=0.75pt,solid,mark=x] 
table[x=no.,y=15,col sep=comma]{KS_effect_no_sigma/res.csv};  
\addplot[color=magenta,line width=0.75pt,solid,mark=x] 
table[x=no.,y=25,col sep=comma]{KS_effect_no_sigma/res.csv};  
\addplot[color=yellow,line width=0.75pt,solid,mark=x] 
table[x=no.,y=30,col sep=comma]{KS_effect_no_sigma/res.csv};  
\end{semilogyaxis}
\end{tikzpicture}

%% file: KS_validation/sensitivities.tex
\begin{tikzpicture}
\begin{axis}[legend entries={$\frac{d\langle \bar{u}\rangle ^{\infty}}{dc}$ \cite{Blonigan2014a},$\frac{d\langle \bar{u} \rangle}{dc}$,
$\frac{d\langle \bar{u^2}\rangle ^{\infty}}{dc}$ \cite{Blonigan2014a},$\frac{d\langle \bar{u^2}\rangle}{dc}$},
 legend pos=north west,  
 legend columns=2,
xlabel=$c$,
   xmin=0,
   xmax=1.2,
	grid=both,
	minor grid style={gray!25},
	major grid style={gray!25},
	width=0.62\textwidth,
	no marks]
\addplot [
    domain=0:1.2, 
    samples=100, 
    color=black,dashed]
{-0.915};
\addplot[color=black,line width=0.75pt,solid,only marks]
table[x=c,y=1,col sep=comma]{KS_validation/sensitivities.csv};  
\addplot [
    domain=0:1.2, 
    samples=100, 
    color=blue,dashed]
{1.792*x+0.042};
\addplot[color=blue,line width=0.75pt,solid,only marks]
table[x=c,y=1,col sep=comma]{KS_validation/sensitivities_u2.csv};  
\addplot[color=black,line width=0.75pt,solid,only marks]
table[x=c,y=2,col sep=comma]{KS_validation/sensitivities.csv};  
\addplot[color=black,line width=0.75pt,solid,only marks]
table[x=c,y=3,col sep=comma]{KS_validation/sensitivities.csv};  
\addplot[color=black,line width=0.75pt,solid,only marks]
table[x=c,y=4,col sep=comma]{KS_validation/sensitivities.csv};  
\addplot[color=black,line width=0.75pt,solid,only marks]
table[x=c,y=5,col sep=comma]{KS_validation/sensitivities.csv};  
\addplot[color=black,line width=0.75pt,solid,only marks]
table[x=c,y=6,col sep=comma]{KS_validation/sensitivities.csv};  
\addplot[color=black,line width=0.75pt,solid,only marks]
table[x=c,y=7,col sep=comma]{KS_validation/sensitivities.csv};  
\addplot[color=black,line width=0.75pt,solid,only marks]
table[x=c,y=8,col sep=comma]{KS_validation/sensitivities.csv};  
\addplot[color=black,line width=0.75pt,solid,only marks]
table[x=c,y=9,col sep=comma]{KS_validation/sensitivities.csv};  
\addplot[color=black,line width=0.75pt,solid,only marks]
table[x=c,y=10,col sep=comma]{KS_validation/sensitivities.csv};  
\addplot[color=black,line width=0.75pt,solid,only marks]
table[x=c,y=11,col sep=comma]{KS_validation/sensitivities.csv};  
\addplot[color=black,line width=0.75pt,solid,only marks]
table[x=c,y=12,col sep=comma]{KS_validation/sensitivities.csv};  
\addplot[color=black,line width=0.75pt,solid,only marks]
table[x=c,y=13,col sep=comma]{KS_validation/sensitivities.csv};  
\addplot[color=black,line width=0.75pt,solid,only marks]
table[x=c,y=14,col sep=comma]{KS_validation/sensitivities.csv};  
\addplot[color=black,line width=0.75pt,solid,only marks]
table[x=c,y=15,col sep=comma]{KS_validation/sensitivities.csv};  
\addplot[color=black,line width=0.75pt,solid,only marks]
table[x=c,y=16,col sep=comma]{KS_validation/sensitivities.csv};  
\addplot[color=black,line width=0.75pt,solid,only marks]
table[x=c,y=17,col sep=comma]{KS_validation/sensitivities.csv};  
\addplot[color=black,line width=0.75pt,solid,only marks]
table[x=c,y=18,col sep=comma]{KS_validation/sensitivities.csv};  
\addplot[color=black,line width=0.75pt,solid,only marks]
table[x=c,y=19,col sep=comma]{KS_validation/sensitivities.csv};  
\addplot[color=black,line width=0.75pt,solid,only marks]
table[x=c,y=20,col sep=comma]{KS_validation/sensitivities.csv};  
\addplot[color=blue,line width=0.75pt,solid,only marks]
table[x=c,y=2,col sep=comma]{KS_validation/sensitivities_u2.csv};  
\addplot[color=blue,line width=0.75pt,solid,only marks]
table[x=c,y=3,col sep=comma]{KS_validation/sensitivities_u2.csv};  
\addplot[color=blue,line width=0.75pt,solid,only marks]
table[x=c,y=4,col sep=comma]{KS_validation/sensitivities_u2.csv};  
\addplot[color=blue,line width=0.75pt,solid,only marks]
table[x=c,y=5,col sep=comma]{KS_validation/sensitivities_u2.csv};  
\addplot[color=blue,line width=0.75pt,solid,only marks]
table[x=c,y=6,col sep=comma]{KS_validation/sensitivities_u2.csv};  
\addplot[color=blue,line width=0.75pt,solid,only marks]
table[x=c,y=7,col sep=comma]{KS_validation/sensitivities_u2.csv};  
\addplot[color=blue,line width=0.75pt,solid,only marks]
table[x=c,y=8,col sep=comma]{KS_validation/sensitivities_u2.csv};  
\addplot[color=blue,line width=0.75pt,solid,only marks]
table[x=c,y=9,col sep=comma]{KS_validation/sensitivities_u2.csv};  
\addplot[color=blue,line width=0.75pt,solid,only marks]
table[x=c,y=10,col sep=comma]{KS_validation/sensitivities_u2.csv};  
\addplot[color=blue,line width=0.75pt,solid,only marks]
table[x=c,y=11,col sep=comma]{KS_validation/sensitivities_u2.csv};  
\addplot[color=blue,line width=0.75pt,solid,only marks]
table[x=c,y=12,col sep=comma]{KS_validation/sensitivities_u2.csv};  
\addplot[color=blue,line width=0.75pt,solid,only marks]
table[x=c,y=13,col sep=comma]{KS_validation/sensitivities_u2.csv};  
\addplot[color=blue,line width=0.75pt,solid,only marks]
table[x=c,y=14,col sep=comma]{KS_validation/sensitivities_u2.csv};  
\addplot[color=blue,line width=0.75pt,solid,only marks]
table[x=c,y=15,col sep=comma]{KS_validation/sensitivities_u2.csv};  
\end{axis}
\end{tikzpicture}

%% file: ill_cond.tex
In Section (\ref{Schur-PC}), a block diagonal preconditioner was presented that can deflate large singular values. While deflating is essential for accelerating convergence, as shown in Figures \ref{Conv dT2 fig} and \ref{eig no sigma fig}, very small singular values are  present for quasi-hyperbolic systems and cause significant problems in solution accuracy and convergence.

Consider again the sensitivity $\sfrac{d\langle \bar{u} \rangle}{dc}$ shown in Figure (\ref{KS sensitivities fig}). There is a constant bias of about $8\%$ for all simulations (20 random $\textbf{u}_0$ for each value of $c$). This bias has been observed previously for the KS in \cite{Blonigan2014a,Blonigan2018}. We can gain helpful insight by considering the analytical solution to the minimisation problem (\ref{MSS system eqn1}), expressed in terms of the singular values and left and right singular vectors of matrix $A$ as follows:
\begin{equation}\label{svd sum eqn}
\underline{\textbf{v}}_{l}=\sum_{i=1}^{l}\left( \frac{u_i^T\underline{\textbf{b}}} {\sigma_i} \right) v_i
\end{equation}

This expression indicates that the minimal norm solution is a linear combination of the right singular vectors $v_i$. The  coefficients $\sfrac{u_i^T\underline{\textbf{b}}}{\sigma_i}$  are obtained by projecting the right hand side $\underline{\textbf{b}}$ to the left singular vectors $u_i^T$ and dividing by the singular value $\sigma_i$. For $l=NK$, i.e.\ using all singular modes, \eqref{svd sum eqn} yields the same solution obtained by solving (\ref{KKT eqn}) iteratively. However, terminating the summation at a value of $l<NK$ allows us to study the effect of a smaller group of singular modes. We can also use a single value of index $i$ to compute the contribution of an individual singular mode to the solution $\underline{\textbf{v}}$, and therefore to the sensitivity. These properties make the decomposition \eqref{svd sum eqn} a very useful tool. 

\begin{figure}[htb!]
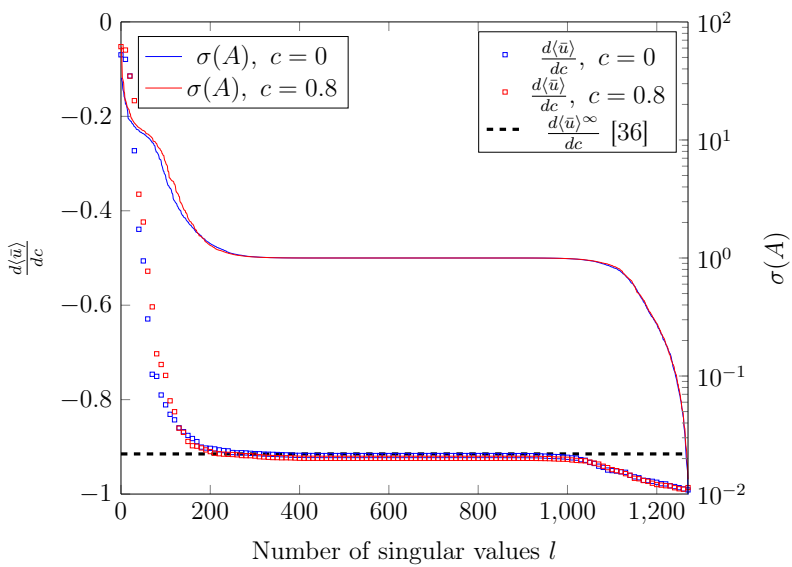

\centering
\includestandalone[width=0.6\textwidth]{KS_validation/sigma_T_100}%
\caption{Sensitivities computed using equation (\ref{svd sum eqn}) for different values of $l$ (KS equation, $T=100$, $N=127$, $K=10$ segments). The solid lines show $\sigma(A)$ (right vertical axis).}
\label{conditioning KS fig}
\end{figure}

Figure (\ref{conditioning KS fig}) shows $\sfrac{d\langle \bar{u} \rangle}{dc}$ for the KS equation, obtained when different values of $l$ are used to terminate the summation \eqref{svd sum eqn}. On the same plot, we superimpose the corresponding singular values (right vertical axis). A very interesting behaviour can be noticed. Summing modes with $\sigma\geq1$ leads to an error in $\sfrac{d \langle \bar{u} \rangle }{dc}$ of less than $1\%$ of $\sfrac{d\langle\bar{u}\rangle^{\infty}}{dc}$. When $l$ is between 250 to 1000,  with $\sigma(A)\approx 1$, the sensitivity remains almost constant. However, including in the summation terms with very small values of $\sigma(A)$ degrades the accuracy of the solution, as seen from the divergence of the squares from the dashed line. The final value is equal to the one shown in Figure (\ref{KS sensitivities fig}). 
\begin{figure}[htb!]
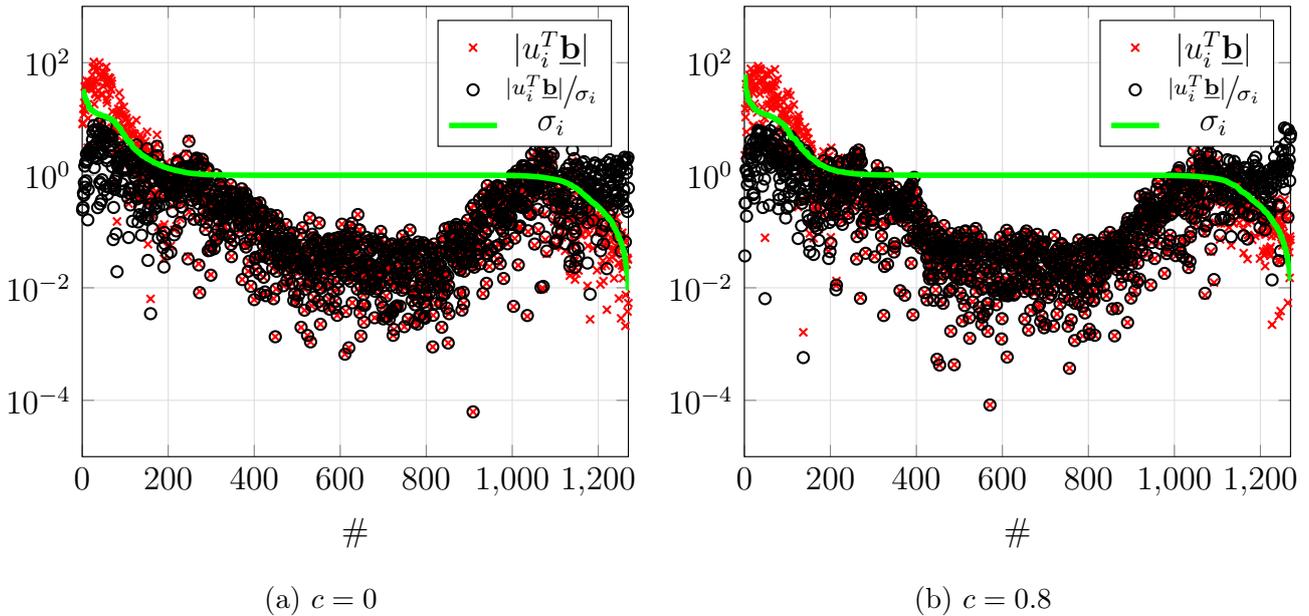

\begin{subfigure}{0.49\textwidth}
 \includestandalone[width=1\textwidth]{KS_validation/coefficients_c0}
 \caption{$c=0$}
 \end{subfigure}
 \begin{subfigure}{0.49\textwidth}
 \includestandalone[width=1\textwidth]{KS_validation/coefficients_c0_8}
  \caption{$c=0.8$}
 \end{subfigure}
  \caption{Spectral coefficients for $T=100$ (KS equation)}
\label{coefficients fig}
\end{figure}

We can gain more insight by plotting the coefficients $\sfrac{|u_i^T\underline{\textbf{b}}|}{\sigma_i}$, the projections $|u_i^T\underline{\textbf{b}}|$, and $\sigma_i$ together. This is known as a discrete Picard plot \cite{ERFANI_2013}, and is shown in Figure (\ref{coefficients fig}). In order to interpret this plot, we recall that very small singular values indicate an ill-conditioned system, i.e.\ the solution is very sensitive to small changes in the right hand side $\underline{\textbf{b}}$. Let's decompose $\underline{\textbf{b}}$ into an unknown error-free part $\underline{\hat{\textbf{b}}}$, and a random error part $\underline{\textbf{e}}$  (for example due to the spatial discretisation and the time advancement scheme), i.e. $\underline{\textbf{b}}=\underline{\hat{\textbf{b}}}+\underline{\textbf{e}}$, with $\|\underline{\textbf{e}}\|_2 \ll \|\underline{\hat{\textbf{b}}}\|_2$. Substituting in \eqref{svd sum eqn} with $l=NK$ we get

\begin{equation}\label{svd sum eqn 2}
\underline{\textbf{v}}=\sum_{i=1}^{NK}\left( \frac{u_i^T \underline{\hat{\textbf{b}}}}{\sigma_i} \right) v_i + \left( \frac{u_i^T\underline{\textbf{e}}} {\sigma_i} \right) v_i = \underline{\hat{\textbf{v}}}+\underline{\textbf{v}_e}
\end{equation}

\noindent This equation indicates that very small $\sigma_i$ can amplify the error component of the solution, $\underline{\textbf{v}_e}=\sum_{i=1}^{NK}\left( \sfrac{u_i^T\underline{\textbf{e}}} {\sigma_i} \right) v_i$. In order to have a meaningful solution, the projection of the error component to the left singular vector, $u_i^T\underline{\textbf{e}}$, should decay to $0$ faster than $\sigma_i$. This is known as the Picard condition \cite{Hansen1994}. Inspection of Figure (\ref{coefficients fig}a) shows that, although the values of $u_i^T\underline{\textbf{b}}$ are quite  spread out, it is clear that they decay by 3-4 orders of magnitude for $i$ between $1-600$. The largest values, of order $O(10^2)$, correspond to small $i$, i.e.\ to the largest singular values. This indicates that the largest contribution to $\underline{\textbf{b}}$ originates from the most rapidly growing modes. For $i>600$, the values of $u_i^T\underline{\textbf{b}}$ remain between $10^{-2}-10^0$, i.e.\ at least 2 orders of magnitude smaller than the maximum. We expect that the small values of $u_i^T\underline{\textbf{b}}$ originate from the random error $\underline{\textbf{e}}$. As long as $\sigma_i \approx 1$, their contribution is innocuous, but when $\sigma_i$ is reduced to $10^{-2}$, they are significantly amplified (as demonstrated from the variation of $\sfrac{u_i^T\underline{\textbf{e}}} {\sigma_i}$) and contaminate the solution. This explains the sensitivity trend in Figure (\ref{KS sensitivities fig}). The same mechanism is valid for $c=0.8$ (panel b), and in fact for all $c$ in the light turbulent regime.  

We have performed simulations for larger $c$, in the convection dominated regime, and we also observed deviations of the computed sensitivity from the solution obtained by finite differences. However, the aforementioned justification could not account for these deviations. For an explanation for the observed differences in that regime, refer to \cite{Blonigan2014a}.

%% file: KS_validation/sigma_T_100.tex
\begin{tikzpicture}
  \begin{axis}[
legend entries={$\frac{d\langle \bar{u}\rangle}{dc} \text{, \,} c=0$,$\frac{d\langle\bar{u}\rangle}{dc} \text{, \,} c=0.8$, $\frac{d\langle\bar{u}\rangle^{\infty}}{dc}$ \cite{Blonigan2014a}},
xlabel=Number of singular values $l$,
ylabel=$\frac{d \langle \bar{u} \rangle }{dc}$,
xmin = 0, xmax = 1270,
ymin = -1, ymax = 0,
axis y line* = left,
	minor grid style={gray!25},
	major grid style={gray!25}, 
    width=0.62\textwidth]
    \addplot[color=blue,only marks,mark=square,mark size=1pt] 
table[x=no.sens,y=sensitivity_c0,col sep=comma]{KS_validation/sigma_T_100.csv};
     \addplot[color=red,only marks,mark=square,mark size=1pt] 
table[x=no.sens,y=sensitivity_c0_8,col sep=comma]{KS_validation/sigma_T_100.csv};
\addplot [
    domain=0:1270,  samples=100, 
    color=black,dashed,line width=1.5pt]
{-0.915};
\end{axis}
   \begin{semilogyaxis}[
      legend entries={$\sigma(A) \text{, \,}  c=0$,$\sigma(A) \text{, \,}  c=0.8$},
      legend pos = north west,
         hide x axis,
         hide y axis,
         axis y line*=right,
          xmin = 0, xmax = 1270,
     ymin = 0.01, ymax = 100,
	minor grid style={gray!25},
	major grid style={gray!25},
    width=0.62\textwidth
]
      \addplot[color=blue,solid,no marks] 
table[x=no.,y=c0_sigma,col sep=comma]{KS_validation/sigma_T_100.csv};
     \addplot[color=red,solid,no marks] 
table[x=no.,y=c0_8_sigma,col sep=comma]{KS_validation/sigma_T_100.csv};
   \end{semilogyaxis}
   \begin{semilogyaxis}[
         xmin=0, xmax=1,
         ymin=0.01, ymax=100,
         legend pos = north west,
         hide x axis,
         axis y line*=right,
         ylabel={$\sigma(A)$},
      minor grid style={gray!25},
	major grid style={gray!25},
    width=0.62\textwidth
]
     \end{semilogyaxis}
\end{tikzpicture}

%% file: KS_validation/coefficients_c0.tex
\begin{tikzpicture}
\begin{semilogyaxis}
[clip mode=individual,
legend entries={$|u_i^T\underline{\textbf{b}}|$,$\sfrac{|u_i^T\underline{\textbf{b}}|}{\sigma_i}$,$\sigma_i$},
  legend pos=north east,
	xlabel=$\#$,
   xmin=0,
   xmax=1270,
   ymax=1e3,
   ymin=1e-5,
	grid=both,
	minor grid style={gray!25},
	major grid style={gray!25}
	]
    \addplot[color=red,line width=0.75pt,only marks,mark=x] 
table[x=no.,y=c0_a,col sep=comma]{KS_validation/sigma_T_100.csv}; 
\addplot[color=black,line width=0.75pt,only marks,mark=o] 
table[x=no.,y=c0_b,col sep=comma]{KS_validation/sigma_T_100.csv}; 
\addplot[color=green,line width=2pt,solid]
table[x=no.,y=c0_sigma,col sep=comma]{KS_validation/sigma_T_100.csv}; 
\end{semilogyaxis}
\end{tikzpicture}

%% file: KS_validation/coefficients_c0_8.tex
\begin{tikzpicture}
\begin{semilogyaxis}
[clip mode=individual,
legend entries={$|u_i^T\underline{\textbf{b}}|$,$\sfrac{|u_i^T\underline{\textbf{b}}|}{\sigma_i}$,$\sigma_i$},
  legend pos=north east,
	xlabel=$\#$,
   xmin=0,
   xmax=1270,
   ymax=1e3,
   ymin=1e-5,
	grid=both,
	minor grid style={gray!25},
	major grid style={gray!25}
	]
    \addplot[color=red,line width=0.75pt,only marks,mark=x] 
table[x=no.,y=c0_8_a,col sep=comma]{KS_validation/sigma_T_100.csv}; 
\addplot[color=black,line width=0.75pt,only marks,mark=o] 
table[x=no.,y=c0_8_b,col sep=comma]{KS_validation/sigma_T_100.csv}; 
\addplot[color=green,line width=2pt,solid]
table[x=no.,y=c0_8_sigma,col sep=comma]{KS_validation/sigma_T_100.csv}; 
\end{semilogyaxis}
\end{tikzpicture}

%% file: Tikhonov.tex
To mitigate the effect of ill-conditioning on the sensitivity $\sfrac{d\bar{J}}{ds}$ for a given system, we can either use (\ref{svd sum eqn}) with a suitable choice of $l$ or regularise the MSS system. Regularisation has a similar filtering effect, i.e. it damps the contribution of very small singular values. The difference is that selecting a particular value of $l$ results in a sharp filter, while regularisation has a smoother filter kernel, as explained in  \cite{Hansen1994}. Tikhonov regularisation is one of the most widely used regularisation techniques. The idea is to solve a regularised version of (\ref{schur eqn}):
\begin{equation}\label{Schur Tikh eqn}
(\gamma I  + S)\underline{\textbf{w}} = \underline{\textbf{b}}
\end{equation}
where $I$ is the identity matrix, and $\gamma>0$ is an appropriately chosen parameter. It has been employed in \cite{Blonigan2018} to improve the conditioning of $S$ for Dowell's Plate equation ($N=4$) and the KS equation ($N\geq 127$). Parameter $\gamma$ shifts all $\mu(S)$ to  $\mu(S)+\gamma$. A large $\gamma$ however relaxes the continuity constraint (\ref{MSS system eqn1}b) and the sensitivities computed through (\ref{MSS sens eqn}) become inaccurate.  If $\gamma$ is chosen adequately, it can improve the accuracy and accelerate the convergence simultaneously (both for very little additional cost). 

In this section, we apply Tikhonov regularization to the preconditioned system (\ref{PC2 sys eqn}), to simultaneously regularize small $\mu(S)$ and to deflate large $\mu(S)$. We have two options: Regularise the original system and then apply preconditioning, i.e. solve
\begin{equation}\label{BDP Tikh eqn_1}
\mathbf{{M}_{BD}}^{(q)}_{(l)}\left( \gamma I+S\right ) \underline{\textbf{w}} =\mathbf{{M}_{BD}}^{(q)}_{(l)}\underline{\textbf{b}} \Longrightarrow 
\left(\gamma \mathbf{{M}_{BD}}^{(q)}_{(l)} +\mathbf{{M}_{BD}}^{(q)}_{(l)} S\right ) \underline{\textbf{w}} =\mathbf{{M}_{BD}}^{(q)}_{(l)}\underline{\textbf{b}}
\end{equation}
or apply preconditioning first and then regularise, i.e. solve 
\begin{equation}\label{BDP Tikh eqn}
\left (\gamma I+\mathbf{{M}_{BD}}^{(q)}_{(l)}S\right ) \underline{\textbf{w}} =\mathbf{{M}_{BD}}^{(q)}_{(l)}\underline{\textbf{b}}
\end{equation}
When the exact preconditioner $M_{(l)}$ (\ref{ext M eqn}) is used and $\gamma$ is much smaller than the largest $l$ retained singular values, the two options are almost identical. For $\mathbf{{M}_{BD}}^{(q)}_{(l)}$, the second option guarantees the clustering of eigenvalues inside the tight range, $\left [ \gamma+\sigma_{min}(\mathbf{{M}_{BD}}^{(q)}_{(l)}S),\gamma+\sigma_{max}(\mathbf{{M}_{BD}}^{(q)}_{(l)}S)\right ]$. Such tight clustering is conducive to rapid convergence of iterative subspace solvers. The second option physically means that we first deflate the rapidly growing modes within each segment and then slightly relax the continuity constraint between segments. This leads to smaller iterations because now the growth in each segment has been suppressed.

\begin{figure}[htb!]
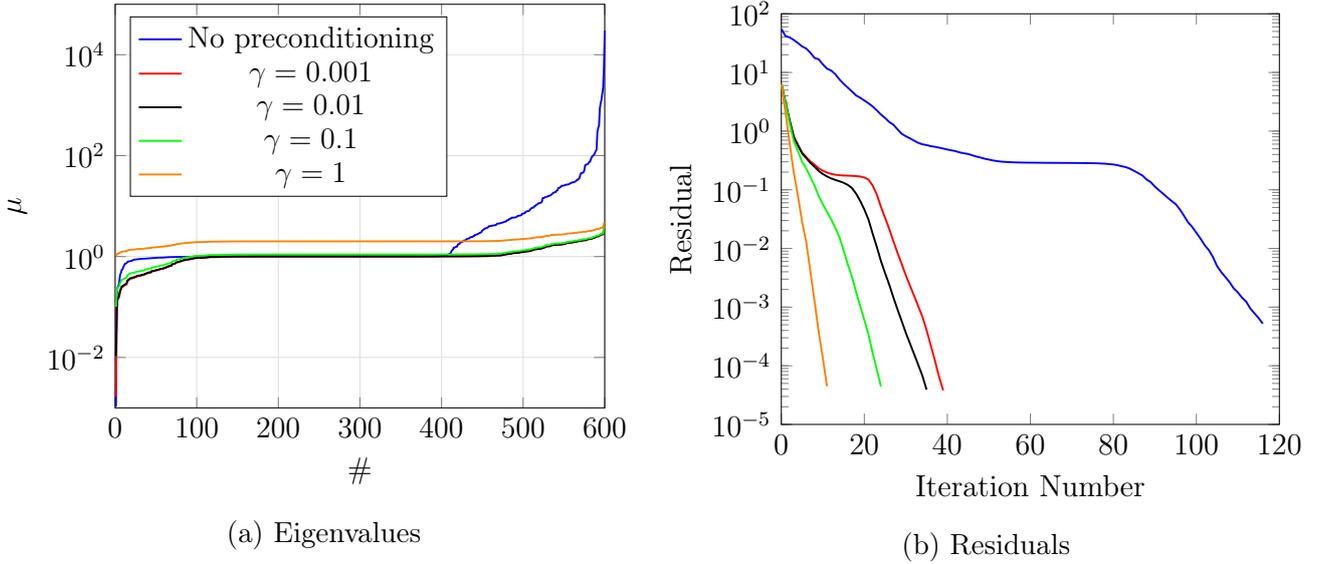

\begin{subfigure}{0.49\textwidth}
\includestandalone[width=0.98\textwidth]{Tikhonov_rho40_T200/eigen_S_T200}
\caption{Eigenvalues}
\end{subfigure}
\begin{subfigure}{0.49\textwidth}
\includestandalone[width=1\textwidth]{Tikhonov_rho40_T200/res_S_T200}
\caption{Residuals}
\end{subfigure}
  \caption{Eigenvalues and residuals of the original system $S$ and the preconditioned system $\gamma I +\mathbf{{M}_{BD}}_{(1)}S$ for different $\gamma$. A Lorenz system trajectory length $T=200$ with $\Delta T=1$ was used for $\rho=40$.}
\label{Tikhonov T200 fig}
\end{figure}

The combined effect of regularization and preconditioning for the Lorenz system can be seen in Figure (\ref{Tikhonov T200 fig}a). Without regularization, the condition number $\kappa(S)\approx 3\times 10^7$ while $\kappa(\mathbf{{M}_{BD}}_{(1)}S)\approx 5.6 \times 10^3$, i.e. a reduction of 4 orders of magnitude in $\kappa$. When using $\gamma=1$, $\kappa(I+\mathbf{{M}_{BD}}_{(1)}S)\approx4$, and convergence is obtained in 12 iterations only (Figure \ref{Tikhonov T200 fig}b). 

\begin{figure}[!htb]
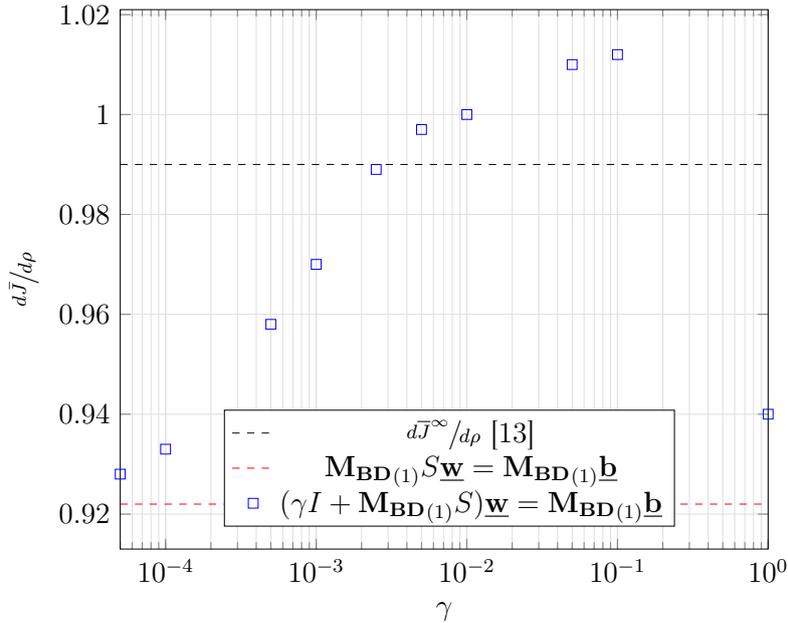

\centering
 \includestandalone[width=0.6\textwidth]{Tikhonov_rho40_T200/sensitivities}
  \caption{Sensitivities for the Lorenz system ($T=200$, $\Delta T=1$ and $\rho=40$) for different $\gamma$.}
\label{Tikhonov T200 1 fig}
\end{figure}

The sensitivities $\sfrac{d\overline{J}}{d\rho}$ computed for different $\gamma$ are shown in Figure (\ref{Tikhonov T200 1 fig}). For reference, $\sfrac{d\overline{J}^\infty}{d\rho}$ \cite{Wang2014} and the solution for $\gamma=0$ (which is affected by the ill-conditioning of $S$) are shown. There is a range $0.001 \leq \gamma \leq  0.1$   which provides a good balance between filtering out the noisy singular values while keeping $\|\underline{\textbf{b}}-A\underline{\textbf{v}}\|_2$ close to zero. In this range, the solution error is $O(\leq 2\%)$. Further increase to $\gamma=1$ relaxes the constraint $A\underline{\textbf{v}}=\underline{\textbf{b}}$ significantly and results in a solution error of $O(5\%)$. A method to estimate the optimal value of $\gamma$ based on the L-curve criterion is proposed in \cite{Calvetti_el_al_1999,Calvetti_el_al_2004}. 

\begin{figure}[h!]
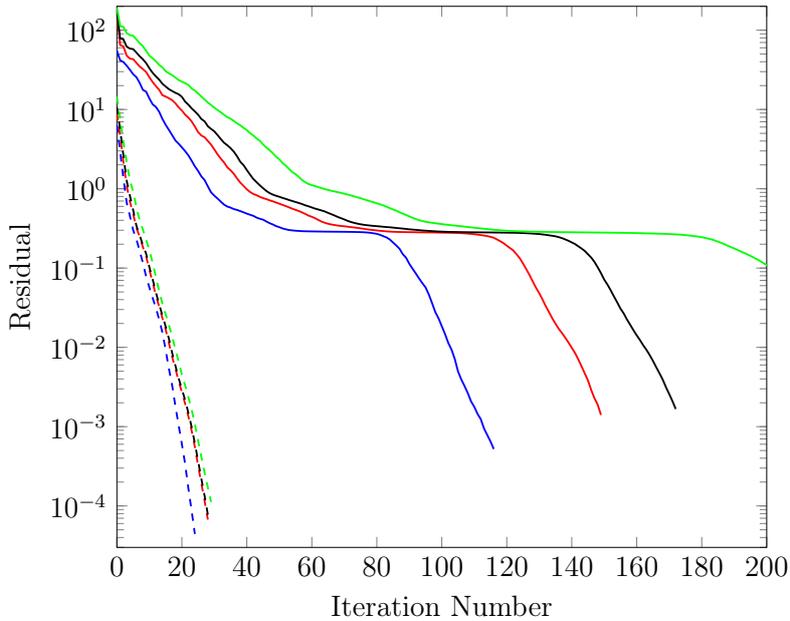

\centering
 \includestandalone[width=0.6\textwidth]{Tikhonov_rho40_dT/res_S_dT}
  \caption{Residuals for $S$ (solid lines) and the BDP system $\gamma I + \mathbf{{M}_{BD}}_{(1)}S$ (dashed lines) with $\gamma=0.1$. The segment size is $\Delta T=1$ and $\rho=40$ (Lorenz system). Blue: $T=200$, red: $T=300$, black: $T=500$, green: $T=1000$.}
\label{Res Tikhonov varyT fig}
\end{figure}

Ideally, the convergence rate should be independent of the number of segments $K$ (and therefore $T$) and the number of degrees of freedom, $N$. This would make MSS applicable to large systems. Figure (\ref{Res Tikhonov varyT fig}) shows that indeed the convergence rate is almost independent of $T$ for the Lorenz system. Theoretical analysis shows that for a linear system with a symmetric, positive definite matrix $\mathcal{A}$, the $\mathcal{A}$-norm of the error at iteration $m$, $\|r_m\|_\mathcal{A}$, satisfies (refer to \cite{Benzi2005})
\begin{equation}\label{conv_rate}
\frac{\|r_m\|_\mathcal{A}}{\|r_0\|_\mathcal{A}} \le 2 \left ( \frac{\sqrt{\kappa \left( \mathcal{A} \right)}-1}{\sqrt{\kappa \left (\mathcal{A} \right )}+1} \right)^m
\end{equation}
This error bound is independent of the number of unknowns, and provided that $\kappa \left(\mathcal{A}\right )$ is suppressed by preconditioning and regularisation, the number of iterations becomes independent of $T$ and $N$. This is clearly demonstrated in Figure (\ref{Res Tikhonov varyT fig}). Tests showed that using equation  \eqref{conv_rate} to predict the number of iterations that will result in a pre-specified normalised residual (set to $10^{-5}$) overestimated the actual iterations needed in practise. This shows that \eqref{conv_rate} indeed provides an upper bound. The bound becomes more accurate as $\kappa \left( \mathcal{A} \right)$ decreases. Table (\ref{table 3}) shows that for the chosen $\gamma=0.1$, $\sfrac{d\overline{J}}{d\rho}\to 0.99$ (the infinite time averaged sensitivity $\sfrac{d\overline{J}^\infty}{d\rho}$). 

\begin{table}[!htb]
\centering
\includestandalone{Tikhonov_rho40_dT/Table}
\caption{A table showing $\sfrac{d\bar{J}}{d\rho}$ for $\rho=40$ (Lorenz system) using different trajectory lengths. A regularization value $\gamma=0.1$ was used.}
\label{table 3}
\end{table}

\begin{figure}[t]
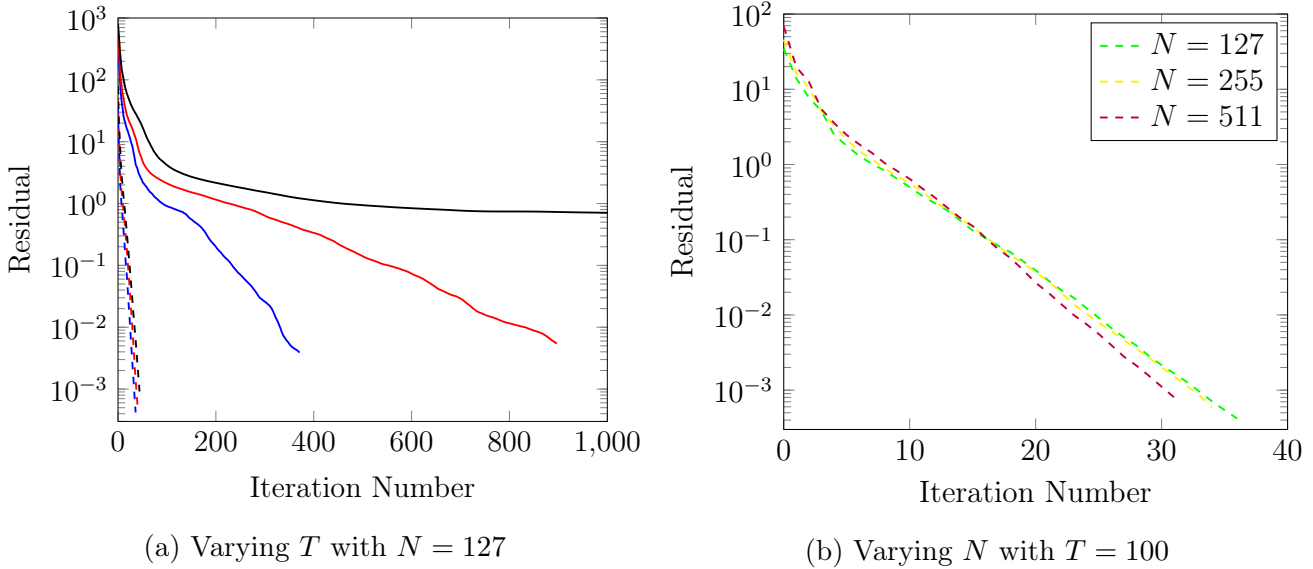

\begin{subfigure}{0.49\textwidth}
 \includestandalone[width=1\textwidth]{KS_Tikhonov/res}
 \caption{Varying $T$ with $N=127$}
\end{subfigure}
\begin{subfigure}{0.49\textwidth}
 \includestandalone[width=1\textwidth]{KS_Tikhonov/res_N}
 \caption{Varying $N$ with $T=100$}
\end{subfigure}
  \caption{Residuals of the original system (solid lines) and of $\gamma I + \mathbf{{M}_{BD}}^{(q)}_{(l)}S$ with $\gamma=0.09$, $q=2$ and $l=15$ (dashed lines). The segment size is $\Delta T=10$ for all cases. Blue: $T=100$, red: $T=200$, black: $T=500$.}
\label{KS Res Tikhonov T fig}
\end{figure}

Figure (\ref{KS Res Tikhonov T fig}) shows the convergence rates for varying $T$ (left panel) and $N$ (right panel) for the Kuramoto Sivashinsky equation. The figure demonstrates again that the combination of regularisation and preconditioning (dashed lines) renders convergence almost independent on $T$ (with $N=127$) and $N$ (with $T=100$). The fast convergence is a direct result of the clustering of eigenvalues in a tight range and the suppression of  $\kappa$.

%% file: Tikhonov_rho40_T200/eigen_S_T200.tex
\begin{tikzpicture}
\begin{semilogyaxis}
[legend entries={No preconditioning,$\gamma = 0.001$,$\gamma = 0.01$,$\gamma = 0.1$,$\gamma = 1$},
  legend pos=north west,
	xlabel=$\#$,
  ylabel=$\mu$,
   xmin=0,
   xmax=600,
   ymax=1e5,
   ymin=1e-3,
	grid=both,
	minor grid style={gray!25},
	major grid style={gray!25},
	no marks]
\addplot[color=blue,line width=0.75pt,solid,mark=x]
table[x=No.,y=eig_nopc,col sep=comma]{Tikhonov_rho40_T200/rho40_gamma0_001.csv};  
\addplot[color=red,line width=0.75pt,solid,mark=x]
table[x=No.,y=eig_pc2,col sep=comma]{Tikhonov_rho40_T200/rho40_gamma0_001.csv};  
\addplot[color=black,line width=0.75pt,solid,mark=x] 
table[x=No.,y=eig_pc2,col sep=comma]{Tikhonov_rho40_T200/rho40_gamma0_01.csv}; 
\addplot[color=green,line width=0.75pt,solid,mark=x]
table[x=No.,y=eig_pc2,col sep=comma]{Tikhonov_rho40_T200/rho40_gamma0_1.csv};  
\addplot[color=orange,line width=0.75pt,solid,mark=x] 
table[x=No.,y=eig_pc2,col sep=comma]{Tikhonov_rho40_T200/rho40_gamma1.csv}; 
\end{semilogyaxis}
\end{tikzpicture}

%% file: Tikhonov_rho40_T200/res_S_T200.tex
\begin{tikzpicture}
\begin{semilogyaxis}
[
	xlabel=Iteration Number,
  ylabel=Residual,
   xmin=0,
   xmax=120,
   ymax=1e2,
   ymin=1e-5,
	no marks]
\addplot[color=blue,line width=0.75pt,solid,mark=x]
table[x=No.res,y=res_nopc,col sep=comma]{Tikhonov_rho40_T200/rho40_gamma0_001.csv};  
\addplot[color=red,line width=0.75pt,solid,mark=x]
table[x=No.res,y=res_pc2,col sep=comma]{Tikhonov_rho40_T200/rho40_gamma0_001.csv};  
\addplot[color=black,line width=0.75pt,solid,mark=x] 
table[x=No.res,y=res_pc2,col sep=comma]{Tikhonov_rho40_T200/rho40_gamma0_01.csv}; 
\addplot[color=green,line width=0.75pt,solid,mark=x]
table[x=No.res,y=res_pc2,col sep=comma]{Tikhonov_rho40_T200/rho40_gamma0_1.csv};  
\addplot[color=orange,line width=0.75pt,solid,mark=x] 
table[x=No.res,y=res_pc2,col sep=comma]{Tikhonov_rho40_T200/rho40_gamma1.csv}; 
\end{semilogyaxis}
\end{tikzpicture}

%% file: Tikhonov_rho40_T200/sensitivities.tex
\begin{tikzpicture}
\begin{semilogxaxis}
[legend style={at={(0.16,0.15)},anchor=west},
legend entries={$\sfrac{d\overline{J}^\infty}{d\rho}$ \cite{Wang2014},$\mathbf{{M}_{BD}}_{(1)}S\underline{\textbf{w}} =\mathbf{{M}_{BD}}_{(1)}\underline{\textbf{b}}$,$(\gamma I+\mathbf{{M}_{BD}}_{(1)}S)\underline{\textbf{w}} =\mathbf{{M}_{BD}}_{(1)}\underline{\textbf{b}}$},
 xlabel=$\gamma$,
  ylabel=$\sfrac{d\bar{J}}{d\rho}$,
   xmin=0.00005,
   xmax=1,
	grid=both,
	minor grid style={gray!25},
	major grid style={gray!25},
    width=0.62\textwidth
]
\addplot [
    domain=0.00005:10, 
    samples=100, 
    color=black,dashed]
{0.99};
\addplot [
    domain=0.00005:10, 
    samples=100, 
    color=red,dashed]
{0.922};
\addplot[color=blue,only marks,mark=square]
table[x=gamma,y=dJ_drho,col sep=comma]{Tikhonov_rho40_T200/sensitivities.csv};  
\end{semilogxaxis}
\end{tikzpicture}

%% file: Tikhonov_rho40_dT/res_S_dT.tex
\begin{tikzpicture}
\begin{semilogyaxis}
[
	xlabel=Iteration Number,
  ylabel=Residual,
   xmin=0,
   xmax=200,
   ymax=2e2,
   ymin=3e-5,
	width=0.62\textwidth,
	no marks]
\addplot[color=blue,line width=0.75pt,solid,mark=x]
table[x=No.res,y=res_nopc,col sep=comma]{Tikhonov_rho40_dT/rho40_T200.csv};  
\addplot[color=blue,line width=0.75pt,dashed,mark=x]
table[x=No.res,y=res_pc2,col sep=comma]{Tikhonov_rho40_dT/rho40_T200.csv};  
\addplot[color=red,line width=0.75pt,solid,mark=x] 
table[x=No.res,y=res_nopc,col sep=comma]{Tikhonov_rho40_dT/rho40_T300.csv}; 
\addplot[color=red,line width=0.75pt,dashed,mark=x]
table[x=No.res,y=res_pc2,col sep=comma]{Tikhonov_rho40_dT/rho40_T300.csv};  
\addplot[color=black,line width=0.75pt,solid,mark=x] 
table[x=No.res,y=res_nopc,col sep=comma]{Tikhonov_rho40_dT/rho40_T500.csv}; 
\addplot[color=black,line width=0.75pt,dashed,mark=x] 
table[x=No.res,y=res_pc2,col sep=comma]{Tikhonov_rho40_dT/rho40_T500.csv}; 
\addplot[color=green,line width=0.75pt,solid,mark=x] 
table[x=No.res,y=res_nopc,col sep=comma]{Tikhonov_rho40_dT/rho40_T1000.csv}; 
\addplot[color=green,line width=0.75pt,dashed,mark=x] 
table[x=No.res,y=res_pc2,col sep=comma]{Tikhonov_rho40_dT/rho40_T1000.csv}; 
\end{semilogyaxis}
\end{tikzpicture}

%% file: Tikhonov_rho40_dT/Table.tex
\begin{tabular}{*5c}
\toprule
 & $T=200$&$T=300$ &  $T=500$ &  $T=1000$   \\ 
\midrule
$\sfrac{d\bar{J}}{ds}$ &1.01&0.97 &0.97&0.99  \\ 
\bottomrule
\end{tabular}

%% file: KS_Tikhonov/res.tex
\begin{tikzpicture}
\begin{semilogyaxis}
[
	xlabel=Iteration Number,
  ylabel=Residual,
   xmin=0,
   xmax=1000,
   ymax=1e3,
   ymin=3e-4,
	no marks]
\addplot[color=blue,line width=0.75pt,solid,mark=x]
table[x=no.,y=no_pc_100,col sep=comma]{KS_Tikhonov/res.csv};  
\addplot[color=blue,line width=0.75pt,dashed,mark=x]
table[x=no.,y=pc_100,col sep=comma]{KS_Tikhonov/res.csv};  
\addplot[color=red,line width=0.75pt,solid,mark=x]
table[x=no.,y=no_pc_200,col sep=comma]{KS_Tikhonov/res.csv};  
\addplot[color=red,line width=0.75pt,dashed,mark=x]
table[x=no.,y=pc_200,col sep=comma]{KS_Tikhonov/res.csv};  
\addplot[color=black,line width=0.75pt,solid,mark=x]
table[x=no.,y=no_pc_500,col sep=comma]{KS_Tikhonov/res.csv};  
\addplot[color=black,line width=0.75pt,dashed,mark=x]
table[x=no.,y=pc_500,col sep=comma]{KS_Tikhonov/res.csv};  
\end{semilogyaxis}
\end{tikzpicture}

%% file: KS_Tikhonov/res_N.tex
\begin{tikzpicture}
\begin{semilogyaxis}
[
legend entries={$N=127$,$N=255$,$N=511$},
	xlabel=Iteration Number,
  ylabel=Residual,
   xmin=0,
   xmax=40,
   ymax=1e2,
   ymin=3e-4,
	no marks]
\addplot[color=green,line width=0.75pt,dashed,mark=x]
table[x=no.,y=N127,col sep=comma]{KS_Tikhonov/res_N.csv};  
\addplot[color=yellow,line width=0.75pt,dashed,mark=x]
table[x=no.,y=N255,col sep=comma]{KS_Tikhonov/res_N.csv};  
\addplot[color=purple,line width=0.75pt,dashed,mark=x]
table[x=no.,y=N511,col sep=comma]{KS_Tikhonov/res_N.csv};  
\end{semilogyaxis}
\end{tikzpicture}

%% file: cost.tex
In this section, we consider the computational cost of the method. This cost is different from the actual wall computing time, as explained later. We use the total number of matrix-vector products involving $\Phi_i$ and $\Phi_i^T$ to quantify the cost. 

The cost of constructing one block of the preconditioner, $\mathbf{M}^{(q)}_{(l),i}$, is $2q(l+2)$, where $l+2$ denotes the selected size of the subspace. The factor $2$ appears because the partial singular value decomposition of $\Phi_i$ requires matrix-vector products with both $\Phi_i$ and $\Phi_i^T$. Since $\mathbf{{M}_{BD}}^{(q)}_{(l)}$ is formed of $K$ blocks, the cost is $2Kq(l+2)$.  One iteration of the subspace solver requires the application of $A$ and $A^T$ once, i.e.\ a total of $2Km$ matrix-vector products are required, where $m$ is the number of iterations. This cost is only approximate, as we have ignored the time taken to orthogonalise the subspace every time a new vector is added.   

Table (\ref{table 4}) shows a cost comparison for different $T$ with and without preconditioning and regularisation (the values correspond to Figure \ref{KS Res Tikhonov T fig}a). We have chosen $q=2$ and $l=15$, so the preconditioner cost is $2Kq(l+2)=68K$. It can be seen that the combination of preconditioning and regularisation results in very significant savings; the cost is reduced by a factor of 35 for $T=500$. Note that the condition number is reduced by between 5 to 7 orders of magnitude (depending on $T$) and remains almost constant. The number of iterations depends very weakly on $T$, and the cost increases almost linearly with $T$ (while keeping $\Delta T$ constant).
\begin{table}[htb!]
\centering
\includestandalone{KS_Tikhonov/table}
\caption{A table showing the cost (number of $\Phi_i$ and $\Phi_i^T$ operations) for the cases shown in Figure (\ref{KS Res Tikhonov T fig}a). The preconditioner was constructed using $q=2$, $l=15$ and $\gamma=0.09$. The relative residual $\sfrac{\|r_m\|_2}{\|r_0\|_2}\approx 1 \times 10^{-5}$ for all cases.}
\label{table 4}
\end{table}
The minimum and maximum eigenvalues are also reported in Table (\ref{table 4}), and this information can be used to assess the individual effects of preconditioning and regularisation. Preconditioning results in a reduction of $\mu_{max}$ by two orders of magnitude (and a corresponding reduction in $\kappa$). Regularisation raises $\mu_{min}$ by three to five orders of magnitude. 

We have not attempted to minimise cost, but there is scope for significant reduction. For example, we have chosen $\gamma=0.09$, which produces a value for the sensitivity $\sfrac{d\langle \bar{u} \rangle^{\infty}}{dc}$ accurate to $1\%$. Increasing to $\gamma=0.25$ and using $q=1$, still gives an acceptable sensitivity (the error is $8\%$), but the number of iterations is reduced to 35, and the total cost is 5800.  

We focused the above analysis on the computational cost. In practical computations, the fully parallel-in-time construction of the preconditioner must be exploited. Each block can be computed independently from all the others. The matrix-vector products involving $A$ and $A^T$ can also be computed in parallel for each segment. Message passing can be overlapped with computations to make the implementation more efficient.

%% file: KS_Tikhonov/table.tex
\centering
\begin{tabular}{*4c}
\toprule
 &  $T=100$ & $T=200$ & $T=500$\\
 &  $(K=10)$ & $(K=20)$ & $(K=50)$\\
\midrule
\multicolumn{4}{c}{Cost of solving$\left(\gamma I+\mathbf{{M}_{BD}}^{(q)}_{(l)}S\right)\underline{\textbf{w}} =\mathbf{{M}_{BD}}^{(q)}_{(l)}\underline{\textbf{b}}$} \\
\hline
$\mu_{max}$, $\mu_{min}$  &$11.87$, $0.090$& $11.87$, $0.090$ &$13.37$, $0.090$\\
Condition number, $\kappa$ &132 & 132&149\\
Number of iterations, $m$ & 38 &42& 46\\
Cost of iterations &760 &1,680 &4,600\\ 
Cost of constructing $\mathbf{{M}_{BD}}^{(q)}_{(l)}$ &680&1,360 &3,400\\
Total cost &1,440 &3,040 &8,000\\ 
\hline
\multicolumn{4}{c}{Cost of solving $S\underline{\textbf{w}} =\underline{\textbf{b}}$} \\
\hline
$\mu_{max}$, $\mu_{min}$  &$3800$, $1.90\times 10^{-4}$ & $3800$, $ 1.80\times 10^{-5}$  &$4900$, $4.90\times 10^{-6}$\\
Condition number, $\kappa$ &$2.0\times 10^7$&$2.1\times 10^8$ &$1.0\times 10^9$\\
Number of iterations, $m$ & 371& 897& 2,790\\
Total cost &7,420 &35,880 &279,000\\ 
\bottomrule
\end{tabular}

%% file: Conclusions.tex
We proposed a block diagonal preconditioner to accelerate the convergence rate for the solution of the linear system arising from the application of the Multiple Shooting Shadowing algorithm. The preconditioner is based on the partial singular value decomposition of the diagonal blocks of the Schur complement. It was applied to the Lorenz system and the Kuramoto Sivashinsky equation. 

The number of singular modes to retain in the partial SVD, $l$, is case dependent. A well-chosen value is required for fast convergence. If the number of positive Lyapunov exponents is known, it can be used to inform the choice of $l$. Strictly speaking however, this is not necessary.  A self-adaptive algorithm is currently being investigated to estimate $l$ that requires no prior knowledge of the number of positive Lyapunov exponents. 

When the preconditioner was combined with a regularisation method, 
the condition number was significantly suppressed, and the convergence rate was found to be weakly dependent on the number of degrees of freedom and the length of the trajectory. The total number of operations was significantly reduced as a result. This paves the way to apply MSS to higher dimensional chaotic systems for sensitivity and optimal control applications. Furthermore, it was shown that a large condition number can affect the accuracy of the computed sensitivity, and can explain an $8\%$ bias in the Kuramoto-Sivashinsky equation in the light turbulence regime.  

Some open issues remain. Apart from the value of $l$ mentioned above, the question on how to choose an appropriate value of the regularisation parameter $\gamma$ that can provide an adequate balance between solution accuracy and rate of convergence is still open. Some ideas to compute the optimal $\gamma$ are provided in \cite{Calvetti_el_al_1999,Calvetti_el_al_2004}. The value is case dependent, and more simulations are required to improve our  understanding, especially for high-dimensional systems. 